\RequirePackage[l2tabu, orthodox]{nag}
\documentclass[12pt]{article}
\usepackage[T1]{fontenc}
\usepackage{microtype}

\usepackage{amsmath}
\usepackage{amssymb}
\usepackage{amsthm}
\usepackage{mathtools}
\usepackage{thmtools}
\usepackage{mathrsfs}

\usepackage{libertine}             

\usepackage{titlesec}

\titleformat*{\section}{\Large\bfseries}

\usepackage[usenames,dvipsnames]{xcolor}
\definecolor{shadecolor}{gray}{0.9}

\usepackage[english]{babel}
\usepackage[parfill]{parskip}
\usepackage{afterpage}
\usepackage{framed}
\usepackage{nicefrac}

{\endMakeFramed}

\newcounter{parcount}

\usepackage{graphicx}
\usepackage[labelfont=bf]{caption}
\usepackage[format=hang]{subcaption}
\usepackage{wrapfig}

\usepackage{booktabs}
\usepackage{multirow}
\usepackage{longtable}

\usepackage[sort]{natbib}
\usepackage{bibunits}

\usepackage[algoruled]{algorithm2e}
\usepackage{listings}
\usepackage{fancyvrb}
\fvset{fontsize=\normalsize}

\usepackage{hyperref}
\hypersetup{colorlinks,linktoc=all}
\usepackage[all]{hypcap}
\hypersetup{citecolor=Violet}
\hypersetup{linkcolor=black}
\hypersetup{urlcolor=MidnightBlue}

\usepackage[nameinlink]{cleveref}

\usepackage[acronym,smallcaps,nowarn]{glossaries}

\usepackage{listings}
\usepackage{lstbayes}
\lstdefinestyle{mystyle}{
    commentstyle=\color{OliveGreen},
    numberstyle=\tiny\color{black!60},
    stringstyle=\color{BrickRed},
    basicstyle=\ttfamily\scriptsize,
    breakatwhitespace=false,
    breaklines=true,
    captionpos=b,
    keepspaces=true,
    numbers=none,
    numbersep=5pt,
    showspaces=false,
    showstringspaces=false,
    showtabs=false,
    tabsize=2
}
\usepackage{enumitem}

\newtheorem{lemma}{Lemma}

\newtheorem{proposition}{Proposition}
\theoremstyle{definition}

\numberwithin{equation}{section}
\newtheorem{thm}{Theorem}[section]

\crefname{lemma}{lemma}{lemmas}
\Crefname{lemma}{Lemma}{Lemmas}
\crefname{thm}{theorem}{theorems}
\Crefname{thm}{Theorem}{Theorems}
\crefname{prop}{proposition}{propositions}
\Crefname{prop}{Proposition}{Propositions}
\crefname{defn}{definition}{definitions}
\Crefname{defn}{Definition}{Definitions}

\crefname{lemma}{lemma}{lemmas}
\Crefname{lemma}{Lemma}{Lemmas}
\crefname{thm}{theorem}{theorems}
\Crefname{thm}{Theorem}{Theorems}
\crefname{prop}{proposition}{propositions}
\Crefname{prop}{Proposition}{Propositions}
\crefname{assumption}{assumption}{assumptions}
\Crefname{assumption}{Assumption}{Assumptions}

\newtheorem{assumption}{Assumption}

\newenvironment{proofsk}{%
  \renewcommand{\proof}{\textbf{Proof sketch.}}\proof}{\endproof}

\newcommand*{\QED}{\hfill\ensuremath{\square}}

\newcommand\independent{\protect\mathpalette{\protect\independenT}{\perp}}
\def\independenT#1#2{\mathrel{\rlap{$#1#2$}\mkern2mu{#1#2}}}

\renewcommand{\mid}{~\vert~}

\newcommand{\cN}{\mathcal{N}}

\newacronym{POMDP}{pomdp}{partially observable {M}arkov decision process}
\newacronym{MDP}{mdp}{{M}arkov decision process}

\newacronym{PNS}{pns}{probability of necessity and sufficiency}
\newacronym{PS}{ps}{probability of sufficiency}
\newacronym{PN}{pn}{probability of necessity}
\newacronym{POC}{poc}{probabilities of causation}
\newacronym{PPCA}{ppca}{probabilistic principal component analysis}

\newacronym{GMM}{gmm}{Gaussian mixture model}

\newacronym{ATT}{att}{average treatment effect on the treated}
\newacronym{ATE}{ate}{average treatment effect}
\newacronym{MAR}{mar}{missing at random}
\newacronym{MSE}{mse}{mean squared error}

\usepackage[margin=1in]{geometry}
\usepackage{tcolorbox}
\usepackage{adjustbox}

\usepackage[defaultlines=3,all]{nowidow}

\title{Large Sample Properties of Matching for Balance}

\author{
  Yixin Wang\\
        UC Berkeley\\
  ywang@eecs.berkeley.edu\\
  \and
  Jos\'{e} R. Zubizarreta\\
  Harvard University\\
  zubizarreta@hcp.med.harvard.edu\\
  }

\date{\today}
\begin{document}
\maketitle
\begin{bibunit}[alp]

\begin{abstract}

Matching methods are widely used for causal inference in observational
studies. Among them, nearest neighbor matching is arguably the most
popular. However, nearest neighbor matching does not generally yield
an average treatment effect estimator that is $\sqrt{n}$-consistent
\citep{abadie2006large}. Are matching methods not
$\sqrt{n}$-consistent in general? In this paper, we study a recent
class of matching methods that use integer programming to directly
target aggregate covariate balance as opposed to finding close
neighbor matches. We show that under suitable conditions these methods
can yield simple estimators that are $\sqrt{n}$-consistent and
asymptotically optimal.

\end{abstract}

Keywords: Causal Inference; Integer Programming; Matching Methods;
Observational Studies; Propensity Score

\section{Introduction}
\label{sec:introduction}

In observational studies, matching methods are widely used for causal
inference. The great appeal of matching methods lies in the
transparency of their covariate adjustments. These adjustments are an
interpolation based on the available data rather than an extrapolation
based on a potentially misspecified model \citep{rubin1973matching,
rosenbaum1989optimal, abadie2006large}. The structure of the data
after matching is also simple (often, a self-weighted sample) so that
statistical inferences and sensitivity analyses are straightforward
\citep{rosenbaum2002observational, rosenbaum2010design1,
rosenbaum2017observation}. Matching methods are commonly used under
the identification assumption of strong ignorability
\citep{rosenbaum1983central} or selection on observables
\citep{imbens2009recent}, but they are also used under the different
assumptions required by instrumental variables (e.g.,
\citealt{baiocchi2010building}) and discontinuity designs (e.g.,
\citealt{keele2015enhancing}).

While there is an extensive literature on matching methods, large
sample characterizations of matching estimators have centered around
nearest neighbor matching only \citep{abadie2006large,
abadie2011bias}. In its simplest form, this algorithm matches each
treated unit to the closest available control in terms of a covariate
distance (e.g., the Mahalanobis distance;
\citealt{rubin1973matching}). 
In an important paper, \citet{abadie2006large} showed that the
resulting difference-in-means estimator is not in general
$\sqrt{n}$-consistent for the average treatment effect when matching with replacement. This estimator
contains a bias that decreases at a rate inversely proportional to the
number of covariates used for matching. As a result, its convergence
can be very slow when matching on many covariates.

Different variants of nearest neighbor matching have been proposed to
address this issue. In one variant, \citet{abadie2011bias} proposed a
class of bias-corrected matching estimators where the missing
potential outcomes are imputed with a regression model. This
imputation corrects the bias of classical nearest neighbor matching.
In another variant, \citet{abadie2016matching} formalized matching on
the estimated propensity score. The estimated propensity score reduces
the matching space into a single dimension.
All of these variants achieve $\sqrt{n}$-consistency. However, in
these cases the faster convergence rate depends on specifying
correctly either the treatment or the outcome model, or restricting the
outcome model to be Lipschitz continuous on the covariates.

\textcolor{black}{Here}, we study a recent class of optimization-based matching
methods that directly target aggregate covariate balance and do not
explicitly model the treatment or the outcome
\citep{zubizarreta2012using,diamond2013genetic,nikolaev2013balance,zubizarreta2014matching}.
These methods formulate the matching exercise as an integer programming problem.
For instance, cardinality matching (\citealt{zubizarreta2014matching}) optimizes the number of matched treated and control units subject
to constraints that approximately balance the empirical distributions
of the covariates. We show that, under suitable conditions, the
resulting difference-in-means treatment effect estimator is
$\sqrt{n}$-consistent, asymptotically Normal, and semiparametrically
efficient. These results show that matching for aggregate covariate
balance can be asymptotically optimal when nearest neighbor matching
with replacement is not.
To our knowledge, this is the first work to show that a matching estimator can be semiparametrically efficient under suitable conditions.

To perform this asymptotic analysis of matching for balance, we
establish a connection between matching and weighting, and view
matching as a form of weighting for covariate balance that both
weights treatment and control units and encodes an assignment between
them. Examples of weighting methods for covariate balance include
\citet{hainmueller2012balancing}, \citet{imai2014covariate},
\citet{zubizarreta2015stable}, \citet{chan2016globally},
\citet{fan2016improving}, \citet{zhao2017entropy},
\citet{athey2018approximate}, \citet{hirshberg2018augmented},
\citet{zhao2019covariate}, and \citet{wang2020minimal}. This
connection between matching and weighting enables us to analyze
matching for balance using asymptotic techniques developed for
weighting.

Despite its connection to weighting, matching for balance retains
some essential features of nearest neighbor matching and other optimal
covariate distance matching algorithms (e.g.,
\citealt{hansen2004full}). In a similar way to distance matching
algorithms, matching for balance can also focus on forming close unit
matches in addition to achieving covariate balance in aggregate. In
fact, matching for balance can be followed by re-matching for
homogeneity in order to not only preserve aggregate covariate balance in the
matched sample, but also minimize total covariate distances between its
matched units (see \citealt{zubizarreta2014matching} for details). As
discussed by \citet{rosenbaum2005heterogeneity} and
\citet{visconti2018handling}, re-matching for homogeneity can improve
the efficiency and sensitivity of certain
matching estimators to unobserved covariates. The main message of this paper is that matching
for balance (along with re-matching for homogeneity) can improve the
large sample properties of the classical difference-in-means estimator
in causal inference, achieving asymptotic optimality under suitable
conditions.

\textcolor{black}{This} paper is organized as follows. In
\Cref{sec:newmatching} we describe the identification assumptions,
matching methods, and the matching estimator. In \Cref{sec:covbaleq}
we present and discuss our main results. 
In \Cref{sec:simulation} we evaluate the empirical performance of the estimator.
In \Cref{sec:conclusion} we
conclude with some remarks. All proofs are shown in the
Supplementary Materials.

\section{Matching for aggregate covariate balance}
\label{sec:newmatching}

In this section, we describe the causal estimation problem and introduce a class of matching methods that target aggregate covariate balance. 
We use the potential outcomes framework for causal inference \citep{neyman1923application,rubin1974estimating}. 
With binary treatments, this framework posits that each unit $i = 1, \ldots, N$ has a pair of potential outcomes $\{ Y_i(0), Y_i(1) \}$, where
$Y_i(1)$ is realized if unit $i$ is assigned to treatment $(Z_i = 1)$
and $Y_i(0)$ is realized if the unit is assigned to control $(Z_i =
0)$. 
\textcolor{black}{Thus, for} each unit $i$, we observe either $Y_i(0)$ or $Y_i(1)$, \textcolor{black}{and the}
observed outcome writes $Y_i = Z_i Y_i(1) + (1-Z_i) Y_i(0)$.
In our setting, the units $i = 1, \ldots, N$ are a random sample from a population of interest and thus the potential outcomes are \textcolor{black}{viewed as} random variables.

Denote $X_i$ as the vector of observed covariates of unit $i$.
These covariates can be continuous or discrete.
Given these
covariates, we assume strong ignorability of the treatment assignment
\citep{rosenbaum1983central}:
$Z_i \independent \{Y_i(0), Y_i(1)\} \mid X_i,$ and
$0 < \Pr(Z_i = 1 \mid X_i) < 1.$
\textcolor{black}{As implied by our notation, we} also require the stable unit
treatment value assumption (SUTVA; \citealt{rubin1980randomization}).

The goal is to estimate the \gls{ATE}, $\mu = \mathbb{E}[Y_i(1) -
Y_i(0)].$ We choose this goal for notational convenience only. 
For example, our
arguments for consistent and efficient estimation of the \gls{ATE} can
be directly extended to the \gls{ATT}, $\textcolor{black}{\mu_t} = \mathbb{E}[Y_i(1) -
Y_i(0)\mid Z_i=1].$

We study matching methods that directly balance the empirical
distributions of the observed covariates. Examples of these methods
are \citet{zubizarreta2012using}, \citet{diamond2013genetic}, and
\citet{nikolaev2013balance}; other related examples \textcolor{black}{include}
\citet{fogarty2015discrete}, \citet{pimentel2015large} and
\citet{kallus2020generalized}. At a high level, these methods aim to
balance the covariates or certain transformations of them that span a
function space (see \citealt{wang2020minimal} for a discussion). We
call these matching methods \emph{matching for balance}. Extending the
formulation in \citet{zubizarreta2014matching}, we study the following matching method
\begin{align}
\text{max.}  \qquad& M\label{eq:matchobjective} \\ 
\text{s.t.} \qquad& m_{ij}\in\{0,1\}, \qquad i,j=1,\ldots,n,\label{eq:matchismatch}\\
\qquad &\sum_{j=1}^n (1-Z_j) m_{ij} = M, \qquad \forall i\in \{i: Z_i = 1\},\label{eq:matchcardcondition}\\
\qquad &\sum_{i=1}^n Z_i m_{ij} = M, \qquad \forall i\in \{j: Z_j = 0\},\label{eq:matchcardcondition_2}\\
\qquad &\sum_{i=1}^n\sum_{j=1}^n Z_iZ_jm_{ij} = \sum_{i=1}^n\sum_{j=1}^n (1-Z_i)(1-Z_j)m_{ij} = 0,\label{eq:matchtrt2ctrl}\\
&\left|\sum_{i=1}^n\sum_{j=1}^n\frac{Z_i(1-Z_j) m_{ij}\{B_k(X_i) - B_k(X_j)\}}{\sum_{i=1}^n\sum_{j=1}^nZ_i(1-Z_j) m_{ij}} \right|<\delta_k, \label{eq:matchcovbalcondition}, k\in[K]\\
&\left|\sum_{i=1}^n\sum_{j=1}^n\frac{(1-Z_i)Z_j m_{ij}\{B_k(X_i) - B_k(X_j)\}}{\sum_{i=1}^n\sum_{j=1}^n(1-Z_i)Z_j  m_{ij}} \right|<\delta_k, k\in[K]\label{eq:matchcovbalcondition_2}
\end{align}
where $m_{ij}$ is a binary decision variable that indicates whether
unit $i$ is matched to unit $j$ (\Cref{eq:matchismatch}).
\textcolor{black}{\Cref{eq:matchcardcondition,eq:matchcardcondition_2}
require each treated unit be matched to $M$ control units, and each
control unit be matched to $M$ treated units, respectively.} \Cref{eq:matchtrt2ctrl}
enforces that each treated unit is not matched to another treated
unit, nor that a control unit is matched to another control unit; in other words, only matches between
different treatment groups are allowed. Finally,
\Cref{eq:matchcovbalcondition,eq:matchcovbalcondition_2} ensure that
the covariate distributions of the matched treated and control units
are balanced. 
In these constraints, the functions $B_k(\cdot)$ are suitable transformations of the covariates. 
Each of them maps the multivariate covariate vector $X_i$ into a suitable summary scalar.
For example, they can be polynomials or wavelets. 
Thus,
\Cref{eq:matchcovbalcondition,eq:matchcovbalcondition_2} constrain the
imbalances in these basis functions in the matched sample up to a
level $\delta_k$. The constant $\delta_k$ is a tuning parameter chosen
by the investigator.
\citet{zhao2019covariate} and \citet{wang2020minimal} describe algorithms to automatically
select the tuning parameter $\delta_k$ in covariate balance
optimization problems like
(\ref{eq:matchobjective})--(\ref{eq:matchcovbalcondition_2}).

\textcolor{black}{As a whole,} optimization problem (\ref{eq:matchobjective})--(\ref{eq:matchcovbalcondition_2}) finds
the largest matched sample with replacement that is balanced according
to the conditions specified in \Cref{eq:matchcovbalcondition} and
\Cref{eq:matchcovbalcondition_2}. 
An interesting feature of this approach is that it accomplishes the task of matching with replacement without predefining the $1:M$ matching ratio, but instead optimizing $M$ from the data at hand subject to the aggregate covariate balance constraints.
We may posit additional
constraints in order to match without replacement: $\sum_{i=1}^n Z_i
m_{ij} \leq 1, \forall j\in \{j: Z_j = 0\}$. In the asymptotic
analyses below, we focus on matching with replacement, but these
analyses can be extended to matching without replacement.

Of course, the above optimization problem of matching for balance may be infeasible. 
For example, this will be the case under practical violations of the positivity assumption\textcolor{black}{; or more specifically, if} there is limited overlap in covariate distributions as characterized by the \textcolor{black}{functions $B_k(\cdot)$.}
In this case, \textcolor{black}{one} is left with the choice to either base the covariate adjustments on a model that \textcolor{black}{goes beyond the support of the data} or to discard some units and possibly change the target of inference \citep{crump2009dealing}.
In this regard, the infeasibility \textcolor{black}{certificate of the matching for balance problem provides} valuable information to characterize the data at hand. 
In the Supplementary Materials, we provide sufficient conditions that guarantee the existence of a solution to the matching for balance optimization problem. 

In order to estimate the \gls{ATE} with matching for balance, we use a simple difference-in-means estimator:
\begin{align}
\hat{\mu}:=\frac{1}{n} & \left[  \sum_{i=1}^nZ_i\left\{Y_i - \frac{\sum_{j=1}^n(1-Z_j)m_{ij}Y_j}{\sum_{j=1}^n(1-Z_j)m_{ij}}\right\} \right.\nonumber\\
&\qquad\qquad\left.+ \sum_{i=1}^n(1-Z_i)\left( \frac{\sum_{j=1}^nZ_jm_{ij}Y_j}{\sum_{j=1}^nZ_jm_{ij}}-Y_i\right)\right].
\label{eq:ATEmatch}
\end{align} 
This estimator computes the average difference between each unit and its matches.
For example, the first term of \Cref{eq:ATEmatch} is the difference
between the outcome of each treated unit $Y_i$ and the mean outcome of
the units it is matched to, $\{Y_j:m_{ij}=1, Z_j=0\}$. Analogously,
the second term is the difference between the outcome of each control unit
$Y_i$ and the mean outcome of its matches,~${\{Y_j:m_{ij}=1, Z_j=1\}}$.

Using \Cref{eq:matchcardcondition} and \Cref{eq:matchcardcondition_2},
we can rewrite this difference-in-means estimator as
\begin{align}
\hat{\mu}
=\frac{1}{n} & \left[\left\{\sum_{i=1}^nZ_iY_i + \sum_{j=1}^n\frac{\sum
_{i=1}^n(1-Z_i)
m_{ij}}{M}Z_jY_j\right\}\right. \nonumber \\
&\qquad\qquad\left.-\left\{\sum_{i=1}^n(1-Z_i)Y_i+\sum_{j=1}^n\frac{\sum
_{i=1}^nZ_im_{ij}}{M}(1-Z_j)Y_j\right\}\right].
\end{align} 
This form implies that each unit $j$ receives weight
$\{1+\sum _{i=1}^n(1-Z_i)m_{ij}/M\}$ if it is treated and weight $\{1+\sum
_{i=1}^nZ_im_{ij}/M\}$ if it is a control. Using
\Cref{eq:matchcardcondition} and \Cref{eq:matchcardcondition_2}, we
can also rewrite the covariate balance constraints in
\Cref{eq:matchcovbalcondition} as
\begin{align*}
\frac{1}{\sum_{i=1}^nZ_i}\left|\sum_{i=1}^nZ_iB_k(X_i) - \sum_{j=1}^n\frac{\sum_{i=1}^nZ_im_{ij}}{M}(1-Z_j)B_k(X_j) \right|&<\delta_k,\nonumber\\
\frac{1}{\sum_{i=1}^n(1-Z_i)}\left|\sum_{i=1}^n(1-Z_i)B_k(X_i) -
\sum_{j=1}^n\frac{\sum
_{i=1}^n(1-Z_i)m_{ij}}{M}Z_jB_k(X_j)
\right|&<\delta_k.
\end{align*}
We observe that the weights of the units in both the constraints and the \gls{ATE}
estimator are functions of the frequencies they are matched,
namely
\begin{align}
\label{eq:weights}
w_T(X_j) &= \frac{\sum
_{i=1}^n(1-Z_i)m_{ij}}{M}  \; \; \; \text{  if $X_j$ is treated,}\\
w_C(X_j) &= \frac{\sum
_{i=1}^nZ_im_{ij}}{M} \; \; \; \text{  if $X_j$ is control.}
\end{align}
We note that the numerator and denominator of the weights $w_T(X_j)$
and $w_C(X_j)$ must be integers due to
\Cref{eq:matchismatch,eq:matchcardcondition}. This restricts the
values that the weights can take. 
\textcolor{black}{Besides this}
constraint, the integer program for matching resembles the convex
optimization problem in covariate balancing weights
\citep{zhao2019covariate,wang2020minimal}.

This connection between matching and weighting allows us to establish
the asymptotic optimality of matching for balance. In the following
section, we show that under suitable conditions, the above
difference-in-means \gls{ATE} estimator is $\sqrt{n}$-consistent,
asymptotically Normal, and semiparametrically efficient. As mentioned,
we focus on the \gls{ATE} for simplicity of exposition. These
consistency and asymptotic normality results readily extend to the
\gls{ATT} because the integer programming problem for the \gls{ATT} is analogous
to that of the \gls{ATE}. The difference is that we match control
units to each treated unit, but not treated units to controls. We note
that calculating the asymptotic variance of the \gls{ATT} estimator
is more nuanced. 
In particular, the semiparametric efficiency bound depends on whether we know the true model for the propensity score.

\section{Asymptotic properties of matching for balance}
\label{sec:covbaleq}

In this section, we show that under standard assumptions, matching for
balance is asymptotically optimal: the resulting \gls{ATE} estimator is $\sqrt{n}$-consistent, asymptotically Normal,
and semiparametrically efficient.

We start by describing the assumptions required. We posit
three sets of conditions on the basis functions  $B(x) = (B_1(x),
\ldots, B_K(x))^\top$, the propensity score function $\pi(x) = \Pr(Z_i=1\,|\,
X_i=x)$, and the mean potential outcome functions
$\mathbb{E}[Y_i(z)\,|\, X_i=x]$ for $z \in \{0,1\}$.

\begin{assumption}
\label{assumption:basis} \sloppy Assume the following conditions on
the basis functions $B(x) = (B_1(x), \ldots, B_K(x))^\top$ hold. There
exist constants $C_0, C_1, C_2 > 0$ such that:
\begin{enumerate}
\item $\sup_{x\in\mathcal{X}}||B(x)||_2\leq C_0K^{1/2}$, where
$\mathcal{X}$ is the domain of the covariates $X$, which is compact.
\item $||\mathbb{E}[B(X_i)^\top B(X_i)]||_2\leq C_1$.
\item $\lambda_{\min}\{\mathbb{E}[B(X_i)B(X_i)^\top]\} > C_2$, where
$\lambda_{\min}\{\mathbb{E}[B(X_i)B(X_i)^\top]\}$ denotes the smallest
eigenvalue of the matrix $\mathbb{E}[B(X_i)B(X_i)^\top].$
\end{enumerate}
\end{assumption}

Assumptions
\ref{assumption:basis}\textcolor{black}{.1}-\ref{assumption:basis}\textcolor{black}{.3}
are standard regularity \textcolor{black}{conditions} on the basis functions. They
restrict their magnitude by the norm of the length-$K$ basis function
vector. These conditions are common in nonparametric sieve estimation
(see Assumption 4.1.6 of \citealt{fan2016improving} and Assumption
2(ii) of \citealt{newey1997convergence}). They are satisfied by many
classes of basis functions including regression spline, trigonometric
polynomial, and wavelet bases \citep{newey1997convergence,
horowitz2004nonparametric, chen2007large, belloni2015some,
fan2016improving}.

\begin{assumption}
\label{assumption:propensity} Assume the following conditions on the
propensity score function $\pi(x) = \Pr(Z_i=1\,|\, X_i=x)$ hold.
\begin{enumerate}
\item There exists a constant $C_3 > 0$ such that $C_3 < \pi(x) <
1-C_3$.
\item There exist \textcolor{black}{vectors} $(\lambda^*_{1T})_{K\times
1}, (\lambda^*_{1C})_{K\times 1}$ such that the true
propensity score function $\pi(\cdot)$ satisfies
$\sup_{x\in\mathcal{X}}|1/\pi(x)-B(x)^\top\lambda^*_{1T}| = O(K^{-r_\pi})$
and $\sup_{x\in\mathcal{X}}|1/\{1-\pi(x)\}-B(x)^\top\lambda^*_{1C}| =
O(K^{-r_\pi})$, where~$r_\pi>1$.
\item There exists a set $\mathcal{M}$ such that the propensity score
function satisfies $1/\pi(x)\in\mathcal{M}$ and
$1/\{1-\pi(x)\}\in\mathcal{M}.$ Moreover, the set $\mathcal{M}$ is a set
of smooth functions such that $\log n_{[]}\{\varepsilon,
\mathcal{M}, L_2(\Pr)\}\leq C_4(1/\varepsilon)^{1/k_1}$, where $n_{[]}\{\varepsilon,
\mathcal{M}, L_2(\Pr)\}$ denotes the covering number of $\mathcal{M}$ by
$\varepsilon$-brackets, $L_2(\Pr)$ is the norm defined as
$||m_1(\cdot)-m_2(\cdot)||_{L_2(\Pr)} =
\mathbb{E}[m_1(X_i) - m_2(X_i)]^2$, $C_4$ is a positive constant,
and~${k_1 > 1/2}$.
\end{enumerate}
\end{assumption}

\Cref{assumption:propensity}\textcolor{black}{.1} requires overlap
between the treatment and control populations. This is part of the
identification assumption described in \Cref{sec:newmatching}.
\Cref{assumption:propensity}\textcolor{black}{.2} is a smoothness
condition on the inverse propensity score function. It requires the
inverse propensity score be uniformly approximable by the basis
functions $B(x) = (B_1(x), \ldots, B_K(x))^\top$. For example, when we
choose the basis functions to be multidimensional splines or power series, this
assumption is satisfied for $r_\pi = s / d$, where $s$ is the number
of continuous derivatives of $1/\pi(x)$ and $d$ is the dimension of
$x$, for $x$ with a compact domain $\mathcal{X}$
\citep{newey1997convergence}.
\Cref{assumption:propensity}\textcolor{black}{.3} constrains the
complexity of the function class to which the inverse propensity score
function belongs. This assumption is satisfied, for example, by the
H\"{o}lder class with smoothness parameter $s$ defined on a bounded
convex subset of $\mathbb{R}^d$ with $s/d > 1$ \citep{van1996weak,
fan2016improving}. This is a key assumption that enables the use of
empirical process techniques in establishing consistency and
asymptotic normality.

\begin{assumption}
\label{assumption:outcome} Assume the following conditions on the mean
potential outcome functions~ $Y_z(x) \stackrel{\Delta}{=}\mathbb{E}[Y_i(z)\,|\, X_i=x]$ for $z \in\{
0,1\}$ \textcolor{black}{hold}.
\begin{enumerate}
\item $\mathbb{E}|Y_i-Y_0(X_i)|<\infty$ and
$\mathbb{E}|Y_i-Y_1(X_i)|<\infty$.

\item $|\mu| < \infty,$ where $\mu =
\mathbb{E}[Y_i(1)-Y_i(0)]$ is the true average treatment effect.

\item There exist $r_y > 1/2$, \textcolor{black}{$(\lambda_{2C}^*)_{K\times
1}$, and $(\lambda_{2T}^*)_{K\times
1}$}
such that $\sup_{x\in\mathcal{X}}|Y_0(x)-B(x)^\top\lambda_{2C}^*| =
O(K^{-r_y})$ and
$\sup_{x\in\mathcal{X}}|Y_1(x)-B(x)^\top\lambda_{2T}^*| =
O(K^{-r_y})$.

\item The potential outcome functions satisfy
$Y_0(\cdot)\in\mathcal{H}$ and $Y_1(\cdot)\in\mathcal{H}$, where
$\mathcal{H}$ is a set of smooth functions satisfying $\log
n_{[]}\{\varepsilon, \mathcal{H}, L_2(\Pr)\}\leq
C_5(1/\varepsilon)^{1/k_2}$, $C_5$ is a positive constant, and $k_2 >
1/2$. As in Assumption \ref{assumption:propensity}\textcolor{black}{.3},
$n_{[]}\{\varepsilon, \mathcal{H}, L_2(\Pr)\}$ denotes the covering number
of $\mathcal{H}$ by $\varepsilon$-brackets, and $L_2(\Pr)$ is the norm
 $||m_1(\cdot)-m_2(\cdot)||_{L_2(\Pr)} = \mathbb{E}[m_1(X_i) -
m_2(X_i)]^2$.

\end{enumerate}
\end{assumption}

Assumptions \ref{assumption:outcome}\textcolor{black}{.1} and
\ref{assumption:outcome}\textcolor{black}{.2} are regularity
conditions on the mean potential outcomes. Assumptions
\ref{assumption:outcome}\textcolor{black}{.3} and
\ref{assumption:outcome}\textcolor{black}{.4} are analogous conditions
to Assumptions \ref{assumption:propensity}\textcolor{black}{.2} and
\ref{assumption:propensity}\textcolor{black}{.3}; they constrain the
smoothness of the mean potential outcome functions and the complexity
of the function class they belong to. 
Under strong ignorability, one may get a rough sense of the function approximation quality (Assumption \ref{assumption:propensity}\textcolor{black}{.3}) by evaluating the prediction error of a fitted outcome model on a holdout sample.
This model explains the observed outcomes in terms of the $K$ basis functions of the observed covariates plus the treatment assignment indicator.
The approximation is likely to be good if the prediction error is small on the holdout data.
That said, such an empirical evaluation will only give us a rough sense of the quality of the approximation, for we only have finite samples.
We note that while no specific modeling assumptions are required for the
inverse propensity score function and the potential outcome functions,
Assumptions \ref{assumption:propensity}\textcolor{black}{.2} and
\ref{assumption:propensity}\textcolor{black}{.3} do require both have
the same form of smoothness, namely that they can all be well approximated
by the same set of basis functions.

\begin{assumption} Assume the following conditions on the matching for
balance problem.
\label{assumption:hyperparam}
\begin{enumerate}

\item $K = o(n^{1/2}).$

\item $||\delta||_2 = O_p[K^{1/2}\{(\log K)/n+K^{-r_\pi}\}]$, where
$\delta = (\delta_1,\ldots, \delta_K).$

\item $n^{\frac{1}{2(r_\pi+r_y-0.5)}} = o(K),$ where $r_\pi, r_y$ are
the smoothness parameters defined in assumptions
\Cref{assumption:propensity}\textcolor{black}{.2} and \Cref{assumption:outcome}\textcolor{black}{.3}.

\end{enumerate}
\end{assumption} 

\Cref{assumption:hyperparam}\textcolor{black}{.1} quantifies the rate at
which the number of basis functions we balance can grow with the
number of units. \Cref{assumption:hyperparam}\textcolor{black}{.2} limits the extent to
which there can be imbalances in the basis functions. Despite these
imbalances, we will show that the optimal large sample properties of
the matching estimator are maintained. \Cref{assumption:hyperparam}\textcolor{black}{.3}
characterizes the growth rates of $K$ and $n$ with respect to the
uniform approximation rates $r_\pi$ and $r_y$.

Now we state the main result of this paper.
\begin{thm} Under
\Cref{assumption:basis,assumption:hyperparam,assumption:outcome,assumption:propensity},
the \gls{ATE} estimator
\begin{align*}
\hat{\mu}&:=\frac{1}{n} \left[  \sum_{i=1}^nZ_i\left\{Y_i - \frac{\sum_{j=1}^m(1-Z_j)m_{ij}Y_j}{\sum_{j=1}^m(1-Z_j)m_{ij}}\right\} + \sum_{i=1}^n(1-Z_i)\left\{ \frac{\sum_{j=1}^mZ_jm_{ij}Y_j}{\sum_{j=1}^mZ_jm_{ij}}-Y_i\right\}\right]
\end{align*} is $\sqrt{n}$-consistent, asymptotically Normal, and
semiparametrically efficient:
\[\sqrt{n}(\hat{\mu} - \mu) \stackrel{d}{\rightarrow} \cN(0,
V_{opt}),\] where $V_{opt}$ equals the semiparametric efficiency bound
\begin{align*} 
V_{opt} = \mathbb{E}\left[\frac{\textrm{Var}[Y_i(1)\mid
X_i]}{\pi(X_i)} + \frac{\textrm{Var}[Y_i(0)\mid
X_i]}{1-\pi(X_i)} +\{\mathbb{E}[Y_i(1)-Y_i(0)\mid
X_i]-\mu\}^2\right]
\end{align*} 
and $\pi(X_i)$ is the propensity score of unit $i$. If in addition
$r_y > 1$ holds, then the estimator
\begin{align*}
\begin{split}
\hat{V}_K = &\sum^n_{i=1} \left[ Z_iw(X_i)Y_i -
\frac{\sum_{i=1}^nZ_iw(X_i)Y_i}{\sum_{i=1}^nZ_i} \right.\\
&- (1-Z_i)w(X_i)Y_i +
\frac{\sum_{i=1}^n(1-Z_i)w(X_i)Y_i}{\sum_{i=1}^n(1-Z_i)} \\
& - \hat{Y}_T(X_i)\times \left\{Z_iw(X_i)
- \frac{1}{\sum_{i=1}^nZ_i}\right\}\\
& \left.+ \hat{Y}_C(X_i)\times \left\{(1-Z_i)w(X_i)
- \frac{1}{\sum_{i=1}^n(1-Z_i)}\right\} \right]^2.
\end{split}
\end{align*} 
is a consistent estimator of the asymptotic variance
$V_{opt}$, where
\begin{align*}
\hat{Y}_T(X_i) &= B(X_i)^\top
\left\{\frac{\sum_{i=1}^nZ_iw(X_i)B(X_i) B(X_i)^\top}{\sum_{i=1}^nZ_i}
\right\}^{-1}
\cdot\left\{ \frac{\sum_{i=1}^nZ_iw(X_i)B(X_i) Y_i}{\sum_{i=1}^nZ_i} \right\}\\
\hat{Y}_C(X_i)  &= B(X_i)^\top \left\{\frac{\sum_{i=1}^n(1-Z_i)w(X_i)B(X_i) B(X_i)^\top}{\sum_{i=1}^n(1-Z_i)}
\right\}^{-1}
\cdot\left\{ \frac{\sum_{i=1}^n(1-Z_i)w(X_i)B(X_i) Y_i}{\sum_{i=1}^n(1-Z_i)} \right\}.
\end{align*}
\label{thm:normality}
\end{thm}

\begin{proofsk}
The proof uses empirical process techniques to analyze \gls{ATE}
estimators as in \citet{fan2016improving} (see also
\citealt{wang2020minimal}). The key challenge in this proof lies in the need to characterize the entire class of
matching solutions in matching for balance. More specifically, the
optimization objective of matching for balance does not involve the
matching solution $m_{ij}$ directly, so it does not correspond to a
unique matching solution. We hence need to study the \gls{ATE}
estimates resulting from all possible matching solutions. In contrast,
the balancing weights \citep{wang2020minimal} and the covariate
balancing propensity score \citep{fan2016improving} both work with
optimization objectives that involve all the weights; these problems
also assume a unique weighting solution with infinite samples.

The proof starts by showing that the implied weights of matching for
balance (Equation (\ref{eq:weights})) approximate the true inverse propensity
score function $\pi(x)^{-1}$. Moreover, this approximation is consistent
due to the balancing constraints
(Equations (\ref{eq:matchcovbalcondition})--(\ref{eq:matchcovbalcondition_2})). The rest
of the proof involves a decomposition of $\hat{\mu}-\mu$
into seven components, where six of them converge to zero in
probability, and the other one is asymptotically Normal and
semiparametrically efficient. Each of the first six components can be
controlled by the bracketing number of the function classes of the
inverse propensity score and the outcome functions.
\Cref{assumption:propensity}\textcolor{black}{.3} and \Cref{assumption:outcome}\textcolor{black}{.4} provide
this control. The full proof is in \Cref{sec:asymptotics_proof} of the
Supplementary Materials. \hfill
\QED
\end{proofsk}

An intuitive explanation of \Cref{thm:normality} relies on two
observations. The first observation is that the \gls{ATE} is an
estimand derived from the entire population rather than from
individual units. 
The asymptotic optimality of \textcolor{black}{our} ATE estimator depends primarily on whether the covariate distribution of the treated units is close in aggregate to that of the control units.
\textcolor{black}{For this type of estimator, how the individual units are matched to each other plays a secondary role.} 
\textcolor{black}{More specifically, the ATE estimator only depends on the number of times each treated (or control) unit is matched.}
For this reason, aggregate covariate
balance is sufficient
for the asymptotic optimality of matching estimators for the
\gls{ATE}. Matching for balance precisely targets this aggregate
covariate balance.
\Cref{eq:matchcovbalcondition,eq:matchcovbalcondition_2} ensure that
the covariate distributions after matching are balanced in aggregate
for the treated and control units.

The second observation is the connection between matching for balance
and covariate balancing weights \citep{hainmueller2012balancing,
imai2014covariate, zubizarreta2015stable, chan2016globally,
fan2016improving, zhao2017entropy, zhao2019covariate,
wang2020minimal}. Both methods formulate the estimation problem as a
mathematical program under similar covariate balancing constraints as
in \Cref{eq:matchcovbalcondition,eq:matchcovbalcondition_2}. Covariate
balancing weights are shown to be asymptotically optimal. Thus,
if matching for balance admits a solution, its implied weights as in
\Cref{eq:weights} are as good as the covariate balancing weights. 
For this reason, under the conditions required by \Cref{assumption:basis,assumption:propensity,assumption:outcome,assumption:hyperparam}, 
matching for balance can also be asymptotically optimal for the difference-in-means estimators for the \gls{ATE}.
When these assumptions do not hold, the nearest neighbor matching estimator for the \gls{ATE} (which is a competing semiparametric estimator) does not similarly achieve the semiparametric efficiency bound.

We conclude this section with a discussion of \Cref{thm:normality} and
its assumptions. Unlike other matching methods that assume a correct
propensity or outcome model, \Cref{thm:normality} studies matching for
balance that posits explicit conditions on covariate balance. Such
conditions are practically appealing because covariate balance is what
is typically checked in practice. Other regularity conditions and
smoothness conditions are standard in nonparametric sieve estimation.

Under these conditions, \Cref{thm:normality} establishes the
asymptotic optimality of the simple difference-in-means estimator for
the \gls{ATE} after matching for aggregate covariate balance.
In practice, these conditions indicate (i) using basis functions such
as power series or wavelets (\Cref{assumption:basis}), (ii)
considering settings with a smooth inverse propensity score and
potential outcome functions, with more continuous derivatives than the
number of covariates
(\Cref{assumption:propensity,assumption:outcome}), and (iii)
restricting the number of balancing basis functions $K$ and the
approximate balance tolerance $\delta$ to scale appropriately with the
number of samples (roughly, $K = O(n^{1/2-\epsilon})$ and
$\delta=O[K^{1/2}\{(\log K)/n + K^{-r_\pi}\}]$, where $\epsilon > 0$ is
a small number; \Cref{assumption:hyperparam}). 
While \Cref{assumption:hyperparam} of \Cref{thm:normality} states that one may balance up to $K = O(n^{1/2-\epsilon})$ basis functions as the sample size goes to infinity, in practice, one should proceed with caution in any given finite sample.
The matching for balance optimization problem may not admit a solution.
Even if a solution exists, the finite sample performance of the resulting estimator may not be ideal.
As \cite{robins1997toward} suggested, nonparametric estimators may suffer from the curse of dimensionality. 
There may not exist a nonparametric estimator with good finite sample performance without the knowledge of the true propensity score.

We note that \citet{abadie2011bias} also devise a matching estimator
that is consistent at the $\sqrt{n}$-rate, but matching for balance
achieves the $\sqrt{n}$-rate in a different way.
\citet{abadie2011bias} correct the bias in nearest neighbor matching
by positing a consistent regression model for the mean potential
outcome function. In contrast, matching for balance avoids this bias
by directly balancing the observed covariates in aggregate. Balancing
covariates in aggregate has been shown to be equivalent to
nonparametric estimation of the inverse propensity score and mean
potential outcome functions
\citep{fan2016improving,zhao2017entropy,hirshberg2018augmented,zhao2019covariate,
wang2020minimal}. This nonparametric approach relieves us from
positing a model for the mean potential outcome function that needs to
be correctly specified. While \Cref{thm:normality} requires certain
conditions on both the propensity score and the potential outcomes
functions, it shows that the matching for balance estimator can
achieve semiparametric efficiency, beyond $\sqrt{n}$-consistency.

Finally, \cite{abadie2012martingale} provide a martingale representation of a widespread nearest neighbor matching estimator and derive its asymptotic distribution.
They decompose the estimator into a martingale term and a conditional bias term. 
Both their and our analysis require the conditional bias term to vanish in order to achieve asymptotic consistency. 
Specifically, \cite{abadie2012martingale} posit regularity conditions under which the conditional bias term converges in probability to zero.
\Cref{thm:normality} utilizes the covariate balance conditions in \Cref{eq:matchcovbalcondition,eq:matchcovbalcondition_2} to ensure that the conditional bias term vanishes.


\section{Simulation studies}
\label{sec:simulation}

Here, we illustrate the empirical performance of matching for balance.
Our simulation study is based on a real dataset about the importance of market access for economic development \citep{redding2008costs}. 
The covariates in this dataset are
non-Gaussian and cannot be characterized by their first two moments.
We focus on a setting where both the propensity score and the outcome
are nonlinear functions of the covariates, and study the \gls{MSE} and coverage probabilities of matching for balance.

To generate the data, we take the actual covariate values from \citet{redding2008costs} and simulate the treatment and \textcolor{black}{outcome values} as follows. 
To simulate the treatment assignment indicator $Z_i$, we first fit a logistic regression model to the original indicator in the dataset. 
Specifically, we fit the model $\Pr(Z_i = 1 \,|\, X_i) = \mathrm{sigmoid}\{(\alpha + \sum_{p=1}^P \beta_j X_{ip} + \sum_{p, p'=1}^P \beta_{pp'} X_{ip}X_{ip'})\}$, where $X_{ip}$ denotes the $p$th observed covariate of unit $i$.
We zero out the estimated coefficients with $p$-values smaller than 0.25 and retain the rest of the coefficients.
We finally generate the treatment assignment indicator $Z_i$ for the simulated dataset via a thresholding model $Z_i = \mathbf{1}{ \{ Z^*_i > 0 \} }$, where $Z^*_i = (\alpha + \sum_{p=1}^P \beta_j X_{ip} + \sum_{p, p'=1}^P \beta_{pp'} X_{ip}X_{ip'}) / c + \text{Unif}(-0.5, 0.5)$, setting $c=50$ to induce limited overlap.

Next, we simulate the potential outcomes $\{Y_i(0), Y_i(1)\}$. Again,
we begin by fitting a linear regression model with all possible
second-order interaction terms to the original treated and control
outcomes in the sample. The model is $Y_i = \alpha' + \sum_{p=1}^P
\beta'_j X_{ip} + \sum_{p, p'=1}^P \beta'_{pp'} X_{ip}X_{ip'} +
\beta'_t T_i + \epsilon_i$, $\epsilon_i \sim N(0,1)$. 
As with the treatment assignment, we zero out all the estimated coefficients with
$p$-values smaller than 0.25, and predict the potential treated and
control outcomes on the entire sample using the fitted model. We then
generate the observed outcomes by $Y_i = Z_iY_i(1)+(1-Z_i)Y_i(0)$.
This data generating process yields a simulated dataset with the same size (122 units) as the original dataset in \citet{redding2008costs} and outcomes determined by three covariates and four associated relevant second-order terms. 
Thus, if we balance up to the $k$th univariate moment of the covariates, we need to balance
$3\cdot k$ basis functions.

With these simulated datasets, we evaluate the \gls{MSE} of matching for balance in estimating the \gls{ATE}. 
We vary the number of basis
functions that we balance by setting the bases to be the moments of
the covariates and increasing their order. We balance the first,
second, and third moments of the covariates. We set the level of
imbalance to be 0.1 standard deviations of the corresponding moment.

\begin{figure}[h!]
\centering
\includegraphics[width=0.5\textwidth]{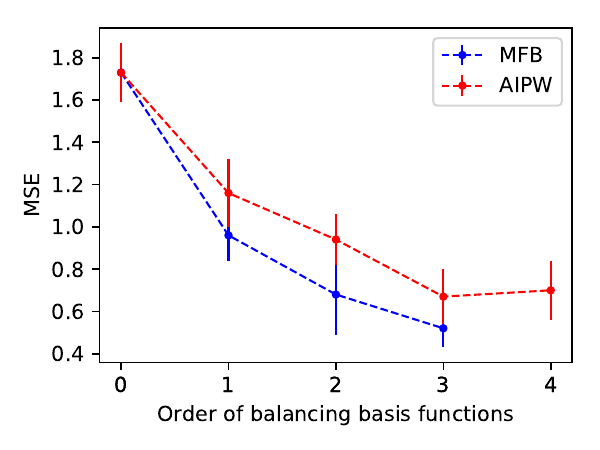}
\caption{The matching for balance (MFB) estimator achieves lower
\gls{MSE} than the augmented inverse probability weighting (AIPW)
estimator when the data has limited overlap. Increasing the number of
balancing basis functions improves the
\gls{MSE} of MFB until its optimization problem
becomes infeasible. The bars indicate $\pm 1$ standard deviation
across 100 simulations. \label{fig:matching-mse}}
\end{figure}

\Cref{fig:matching-mse} shows that as we increase the number of
balancing basis functions, the MSE of matching for balance decreases.
\textcolor{black}{Also, because the simulated dataset has limited
overlap, matching for balance achieves lower MSE than the standard
augmented inverse propensity weighted estimator (AIPW;
\citealt{robins1994estimation}) as the approximate balance constraints
trade variance for bias.} 
Nonetheless, when covariate distributions have limited overlap, this improvement in \gls{MSE} comes at a cost.
As we show below, the resulting confidence intervals may exhibit lower than nominal coverage due to approximate covariate balance.
This suggests exploring separate imbalance tolerances for estimation and for inference.
Finally, because of the nonlinearity of the inverse propensity score, including third order basis functions improves \gls{MSE} despite the data generating process only involves second order terms.

\Cref{fig:matching-mse} also corroborates the
discussion in \Cref{sec:covbaleq} that matching for balance may not
admit a solution if we aim to balance too many basis functions. For
example, in \Cref{fig:matching-mse}, matching for balance admits a
solution if we balance the first three moments of the covariates, but
not if we also try to balance the fourth moment for the given tolerance level of
covariate imbalance. 
We note that in such cases, $\sqrt{n}-$consistency may not hold because
the asymptotic results derived in \Cref{sec:covbaleq} are only applicable when
$K=o(n^{1/2})$. These asymptotic results do not apply if matching for
balance does not admit \textcolor{black}{a solution} when $K$ increases in this order and \textcolor{black}{not enough matches can be found.}

\begin{figure}[h!]
\centering
\begin{subfigure}[b]{0.45\textwidth}
\includegraphics[width=\textwidth]{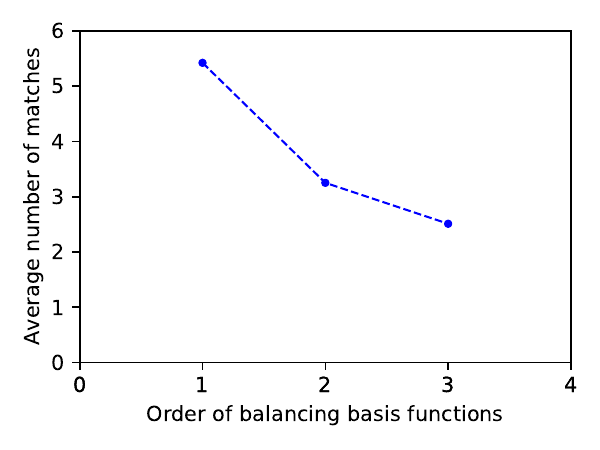}
\caption{Average number of matches\label{fig:nummatch}.}
\end{subfigure}
\begin{subfigure}[b]{0.45\textwidth}
\includegraphics[width=\textwidth]{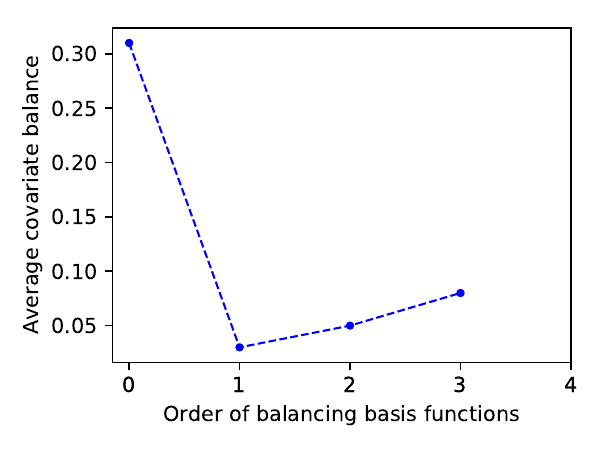}
\caption{Average covariate balance\label{fig:covbal}.}
\end{subfigure}
\caption{Average number of matches and average standardized
covariate balance of matching for balance. While \textcolor{black}{both measures decrease as we balance more basis functions},
the \gls{MSE} still improves because the covariate distributions of
the treatment and control groups are more similar in ways that are
relevant to the propensity score and outcome models.
\label{fig:matching-outcome}}
\end{figure}

\Cref{fig:nummatch} shows that as we increase the number of balancing basis functions, the average number of matches $M$ decreases. 
The reason is that balancing more basis functions implies solving a more constrained optimization problem; hence,  the average $M$ decreases. 
A similar phenomenon \textcolor{black}{appears} in \Cref{fig:covbal} \textcolor{black}{where} the average standardized covariate balance
\[\frac{1}{K}\sum_{k=1}^K\left[
\left|\sum_{i=1}^n\sum_{p=1}^n\frac{Z_i(1-Z_j) m_{ij}\{B_k(X_i) -
B_k(X_j)\}}{\sum_{i=1}^n\sum_{p=1}^nZ_i(1-Z_j) m_{ij}} \right| /
\mathrm{sd}\{B_k(X_i)\}\right] \] 
increases with the number of balancing basis functions. (In
\Cref{fig:covbal}, the order of balancing basis functions equal to
zero represents covariate balance \textcolor{black}{before matching}.) 
While both metrics (the average number of matches and the average absolute standardized mean difference in covariates) \textcolor{black}{decrease} as we increase the number of basis functions that we balance,
the \gls{MSE} still improves because the treatment and control groups
after matching are more similar in ways that are relevant to the propensity score and outcome models. 
This illustrates the importance of balancing more basis functions than just means when the propensity score and outcome models \textcolor{black}{are non-linear} on the covariates. 
While
balancing a high number of basis functions can be difficult with most
matching methods (as they do not directly target covariate balance),
with matching for balance, covariate balance on the basis functions is
obtained by \textcolor{black}{construction}. 
Subject to these covariate balance requirements,
the matching ratio $M$ in matching for balance is optimized, resulting
in the largest possible $1:M/M:1$ matching ratio for the data at hand.

\begin{figure}[htb]
\centering
\includegraphics[width=0.5\textwidth]{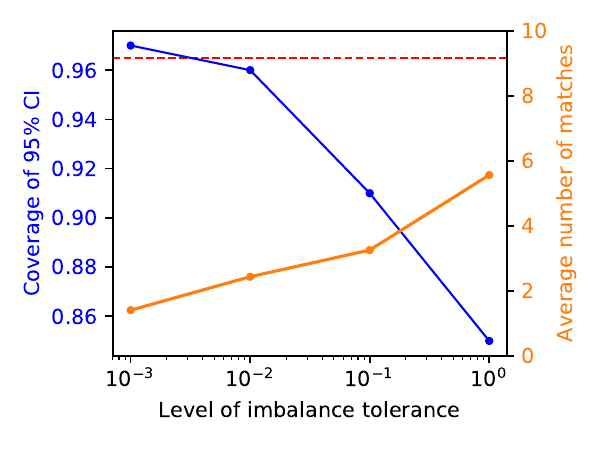}
\caption{The confidence intervals of matching for balance achieve
nominal coverage when the level of imbalance tolerance
\textcolor{black}{$\delta$ is small. Increasing $\delta$} trades
variance for bias but can degrade the coverage probabilities of the
corresponding confidence intervals. In contrast, AIPW (the dashed red
line) achieves close-to-nomimal coverage.
\label{fig:matching-coverage}}
\end{figure}

Finally, we \textcolor{black}{evaluate the coverage} probabilities of
matching for balance and AIPW. \textcolor{black}{We focus} on
balancing the first two moments of the covariates
\textcolor{black}{and vary the imbalance} tolerance $\delta$ such that
$\delta_k = \delta \cdot \mathrm{sd} \{ B_k(X_i) \}$.
\textcolor{black}{In \Cref{fig:matching-coverage}, we show the
coverage probabilities of the 95\% confidence intervals constructed
based on \Cref{thm:normality}.} \textcolor{black}{The figure shows}
that when the imbalance tolerance is small ($\delta \leq 0.01$), the
confidence intervals \textcolor{black}{have} close to nominal
coverage. As the imbalance tolerance increases, the average number of
matches increases but the coverage probability degrades as the matched
sample exhibits worse balance. In contrast, the AIPW estimator achieves close to nominal coverage. In such cases, matching for balance trades \textcolor{black}{the imbalance tolerance $\delta$ for the matching ratio $M$ in order to exchange bias for variance and achieve lower MSE, but it can compromise coverage probabilities.}

An interesting direction for future research is to use a larger value of $\delta$ for estimation, and a smaller value for inference, such that confidence intervals are not necessarily centered at the point estimate.
This would be analogous to what is sometimes done in the analysis of regression discontinuity designs \citep{calonico2014robust}. 
In general, how to select the parameter $\delta$ is an open question in causal inference. 
In the context of weighting, two proposals are provided by \cite{zhao2019covariate} and \cite{wang2020minimal}, where $\delta$ is selected to generalize covariate balance across cross-validation and bootstrap samples, respectively.


\section{Concluding Remarks}
\label{sec:conclusion}

\textcolor{black}{We have analyzed} a recent class of matching methods that \textcolor{black}{targets aggregate covariate balance.
After all, covariate balance} is the main diagnostic that investigators carry out in practice. 
As discussed, matching for balance does not preclude also finding close unit matches\textcolor{black}{, because} matching for balance can be followed by matching for homogeneity in order to minimize covariate distances between matched units, while preserving aggregate covariate balance (see \citealt{zubizarreta2014matching}).

Under \textcolor{black}{suitable} conditions, we \textcolor{black}{have} shown that this class of matching methods yields a
simple difference-in-means estimator that is asymptotically optimal: it is $\sqrt{n}$-consistent, asymptotically normal, and semiparametrically efficient.
\textcolor{black}{As discussed in the simulation study section, these conditions can be stringent in practice because they require the imbalance tolerance $\delta$ to decrease as the sample size increases and there needs to exist a matching solution for such values of $\delta$.}
These results complement the
fundamental results by \cite{abadie2006large}, who showed that a
similar estimator is not in general $\sqrt{n}$-consistent for nearest
neighbor matching with replacement.

Matching for balance exemplifies how tools from modern optimization (e.g., \citealt{junger200950} and \citealt{bixby2012brief}) can play a central role in the design of observational studies in general (e.g., \citealt{rosenbaum2002observational} and \citealt{imbens2015causal}) and in matching for covariate balance in particular (e.g., \citealt{zubizarreta2013stronger} and \citealt{keele2015enhancing}). 
A natural future research direction is to augment matching for balance as in doubly robust estimators (\citealt{robins1994estimation}; see also \citealt{rubin1979using,abadie2011bias,athey2018approximate} for related approaches). 
Another promising direction is to build on the methods described in
\cite{rosenbaum2017observation} on evidence factors and
sensitivity analysis for interpretable analyses of matched
observational studies.
\clearpage
{\putbib[mybibliography17]}
\end{bibunit}

\clearpage
\begin{bibunit}[alp]
\appendix

\allowdisplaybreaks
\clearpage
{\Large\textbf{Supplementary Materials}}

\section{Proof of Theorem 3.1}
\label{sec:asymptotics_proof}

We first revisit the notation. We consider $n$ units randomly drawn
from a population. For each unit $i$, we observe its treatment
indicator $Z_i$, its realized outcome $Y_i$, and its pretreatment
covariates $X_i$. Given the covariate vectors $X_i, i =  1, \ldots,
n$, we consider a vector of $K$ functions of the covariate vector
$B(x) = (B_1(x),\ldots, B_K(x))$.

We are interested in estimating the \gls{ATE} $\mu =
\mathbb{E}[Y_i(1) - Y_i(0)].$ Let $m_{ij}$ be a binary indicator of
whether unit $i$ is matched with unit $j$, $i,j=1, \ldots, n$. The
general form of the matching for
\gls{ATE} is
\begin{align}
\text{max.}  \qquad& M\label{eq:suppmatchobjective} \\ 
\text{s.t.} \qquad& m_{ij}\in\{0,1\}, \qquad i,j=1,\ldots,n,\label{eq:suppmatchismatch}\\
\qquad &\sum_{j=1}^n (1-Z_j) m_{ij} = M, \qquad \forall i\in \{i: Z_i = 1\},\label{eq:suppmatchcardcondition}\\
\qquad &\sum_{i=1}^n Z_i m_{ij} = M, \qquad \forall i\in \{j: Z_j = 0\},\label{eq:suppmatchcardcondition_2}\\
\qquad &\sum_{i=1}^n\sum_{j=1}^n Z_iZ_jm_{ij} = \sum_{i=1}^n\sum_{j=1}^n (1-Z_i)(1-Z_j)m_{ij} = 0,\label{eq:suppmatchtrt2ctrl}\\
&\left|\sum_{i=1}^n\sum_{j=1}^n\frac{Z_i(1-Z_j) m_{ij}[B_k(X_i) - B_k(X_j)]}{\sum_{i=1}^n\sum_{j=1}^nZ_i(1-Z_j) m_{ij}} \right|<\delta_k, \qquad k = 1, \ldots, K, \label{eq:suppmatchcovbalcondition}\\
&\left|\sum_{i=1}^n\sum_{j=1}^n\frac{(1-Z_i)Z_j m_{ij}[B_k(X_i) - B_k(X_j)]}{\sum_{i=1}^n\sum_{j=1}^n(1-Z_i)Z_j  m_{ij}} \right|<\delta_k, \qquad k = 1, \ldots, K, \label{eq:suppmatchcovbalcondition_2}
\end{align}

To estimate the \gls{ATE}, we compute the difference between the mean
outcome of the treated units and that of the matched control units:
\begin{align}
\hat{\mu}&:=\frac{1}{n} \left[  \sum_{i=1}^nZ_i\left(Y_i - \frac{\sum_{j=1}^m(1-Z_j)m_{ij}Y_j}{\sum_{j=1}^m(1-Z_j)m_{ij}}\right) \right.\nonumber\\
&\left.\qquad\qquad+ \sum_{i=1}^n(1-Z_i)\left( \frac{\sum_{j=1}^mZ_jm_{ij}Y_j}{\sum_{j=1}^mZ_jm_{ij}}-Y_i\right)\right].
\label{eq:suppatematch}
\end{align} 

We will give conditions such that the \gls{ATE} estimator is
consistent and asymptotically normal and semiparametric efficient. The
asymptotics for the \gls{ATE} estimator from its corresponding
optimization problem can be proved similarly.

The asymptotic properties of $\hat{\mu}$ mainly depends on two of the
conditions in the optimization problem, namely
\Cref{eq:suppmatchcardcondition} and
\Cref{eq:suppmatchcovbalcondition}.

We first rewrite the \gls{ATE} estimator (\Cref{eq:suppatematch}) as a
weighting estimator:
\begin{align}
\hat{\mu}
=&\frac{1}{n} \left[\sum_{i=1}^nZ_iY_i + \sum_{j=1}^n\frac{\sum
_{i=1}^n(1-Z_i)
m_{ij}}{M}Z_jY_j\right.\nonumber\\
&\left.\qquad\qquad-\sum_{i=1}^n(1-Z_i)Y_i-\sum_{j=1}^n\frac{\sum
_{i=1}^nZ_im_{ij}}{M}(1-Z_j)Y_j\right].
\end{align} 

We can similarly rewrite the balancing conditions
(\Cref{eq:suppmatchcovbalcondition} and
\Cref{eq:suppmatchcovbalcondition_2}) as
\begin{align}
\frac{1}{\sum_{i=1}^nZ_i}\left|\sum_{i=1}^nZ_iB_k(X_i) - \sum_{j=1}^n\frac{\sum_{i=1}^nZ_im_{ij}}{M}(1-Z_j)B_k(X_j) \right|&<\delta_k,\\
\frac{1}{\sum_{i=1}^n(1-Z_i)}\left|\sum_{i=1}^n(1-Z_i)B_k(X_i) -
\sum_{j=1}^n\frac{\sum
_{i=1}^n(1-Z_i)m_{ij}}{M}Z_jB_k(X_j)
\right|&<\delta_k.
\end{align}

Rewriting the difference-in-means estimator and the balancing
conditions implies the following weights assigned to each unit. These
implied weights are functions of the covariates,
\begin{align}
w(X_j) = \begin{cases}
\frac{1}{n}+\frac{\sum
_{i=1}^n(1-Z_i)
m_{ij}}{nM} &\text{ if } Z_j=1,\\
\frac{1}{n}+\frac{\sum
_{i=1}^n Z_i
m_{ij}}{nM} &\text{ if } Z_j=0.
\end{cases}
\end{align}
If two subjects $j, j'$ receive the same treatment assignment
$Z_j=Z_{j'}$ and share the same covariates $X_{j'}=X_{j}$ but
different weights $w(X_{j'})\ne w(X_j)$, we define the weights as an
average
\begin{align}
w(X_j) = \frac{1}{2}(w(X_j) + w(X_{j'})).
\end{align}
Formulating a weight function this way does not affect the value of
either $\hat{\mu}$ or the balancing conditions.

For technical convenience, we also assume that $w(x)$ is a piecewise
constant function in $x$ with jumps of size at most $O_p(K((\log K)/n +
K^{-r_{\pi}})/n)$. More technically, there exists $\delta_0>0$ such that
\[|Z_iw(X_i) - Z_jw(X_j)|<O_p(K((\log K)/n +
K^{-r_{\pi}})/n)\] and
\[|(1-Z_i)w(X_i) -(1-Z_j)w(X_j)|<O_p(K((\log K)/n +
K^{-r_{\pi}})/n)\] for
any $|X_i-X_j|\leq
\delta_0$. This assumption roughly implies that
units with similar covariates are matched for a similar number of
times. Under this assumption, the function $w(x)$ can be uniformly
approximated by an arbitrarily smooth function up to error
$O_p(K((\log K)/n +
K^{-r_{\pi}})/n)$ due to the compactness of
$\mathcal{X}$.

This assumption is satisfied when $1/M = O_p(K((\log K)/n +
K^{-r_{\pi}})/n)$, for example, under the overlap condition posited by
\Cref{prop:solutionexist}. If so, it suffices to guarantee that there
exists $\delta_0>0$ such that $|Z_iw(X_i) - Z_jw(X_j)|<1/M$ and
$|(1-Z_i)w(X_i) - (1-Z_j)w(X_j)|<1/M$ for any
$|X_i-X_j|\leq\delta_0$. In other words, if two units have similar
covariates, then the number of times of each of them being matched can
differ by at most one. We can ensure this condition by considering
weights obtained by first solving a minimization problem
\begin{align*}
\min_{w^0(\cdot)}&\sum_{k=1}^K \left|\left|\frac{\sum_{i=1}^nZ_iB_k(X_i) - \sum_{j=1}^nw^0(X_j)B_k(X_j)}{\delta_k \sum_{i=1}^nZ_i}\right|\right|^2\\
&+\sum_{k=1}^K \left|\left|\frac{\sum_{i=1}^n(1-Z_i)B_k(X_i) - \sum_{j=1}^nw^0(X_j)B_k(X_j)}{\delta_k \sum_{i=1}^n(1-Z_i)}\right|\right|^2, 
\end{align*}
then rounding the solution to the nearest $1/M$ unit, i.e. $w(X_j) =
\frac{[w^0(X_j)\cdot M]}{M}.$ We note that the $w(\cdot)$ is a smooth
function as a solution to the above minimization problem. Under the
conditions of \Cref{prop:solutionexist}, $w(\cdot)$ satisfies the
piece-wise constant condition above.

Rewriting the balancing conditions using the $w(X_j)$ notation, we have
\begin{align}
\frac{1}{\sum_{i=1}^nZ_i}\left|\sum_{i=1}^nZ_iB_k(X_i) - \sum_{j=1}^n(nw(X_j)-1)\cdot(1-Z_j)B_k(X_j) \right|&<\delta_k,\\
\frac{1}{\sum_{i=1}^n(1-Z_i)}\left|\sum_{i=1}^n(1-Z_i)B_k(X_i) -
\sum_{j=1}^n (nw(X_j)-1)\cdot Z_jB_k(X_j)
\right|&<\delta_k.
\end{align}

We first prove the following lemma for the implied weights of the
treated units. An analogous lemma holds for the implied weights of the
control units.

\begin{lemma}
\label{lemma:wtconsistency} Let $\pi(x)$ be the true propensity score
function. Under \Cref{assumption:basis} and
\Cref{assumption:propensity}, we have, for all treated units,
\begin{enumerate}
    \item $\sup_{x\in\mathcal{X}}|nw(x) - 1/\pi(x)| =O_p(K((\log K)/n+K^{-r_\pi}))=o_p(1)$,
    \item $\left\|nw(x) - 1/\pi(x)\right\|_{P,2} =O_p(K((\log K)/n+K^{-r_\pi}))= o_p(1)$.
\end{enumerate}
\end{lemma}

\sloppy
\paragraph{Proof.}
For notation convenience, we define $\mathbf{X} = (X_1, \ldots, X_n)$.
This leads to $B(\mathbf{X}) = (B(X_1), \ldots, B(X_n))^\top$ and
$w(\mathbf{X}) = (w(X_1), \ldots, w(X_n))^\top$. 

We then define $\mathbf{X}_\mathcal{C} = (\{X_i\}_{Z_i=0})^\top$ as
the covariate matrix of the control units. This leads to a similar
basis matrix $B(\mathbf{X}_\mathcal{C})$, that is a $C$ by $K$ matrix,
where $C = \sum_{j=1}^n (1-Z_j)$ is the number of control units, and
$K$ is the number of basis functions we match. Similarly define
$\mathbf{X}_\mathcal{T} = (\{X_i\}_{Z_i=1})^\top$ as the covariate
matrix of the treated units. $B(\mathbf{X}_\mathcal{T})$ as the
covariate basis matrix of the treated units with size $T$ by $K$,
where $T =\sum_{j=1}^n Z_j$ is the number of treated units.

Denote $\lambda^\dagger = [B(\mathbf{X}_\mathcal{T})^\top
B(\mathbf{X}_\mathcal{T})]^{-1}[B(\mathbf{X}_\mathcal{T})^\top
w(\mathbf{X}_\mathcal{T})]$. It is the least square projection of the
weights $w(\mathbf{X}_\mathcal{T})$ onto the space of basis functions
$B(\mathbf{X}_\mathcal{T})$.

Below we prove that $\left\|\lambda^\dagger-\lambda_{1T}^*\right\|_2 =
O_p(K^{1/2}((\log K)/n+K^{-r_\pi})/n)$:
\begin{align*}
&\left\|\lambda^\dagger-\lambda_{1T}^*\right\|_2\\
=&\left\|[B(\mathbf{X}_\mathcal{T})^\top B(\mathbf{X}_\mathcal{T})]^{-1}[B(\mathbf{X}_\mathcal{T})^\top
w(\mathbf{X}_\mathcal{T})] - \lambda_{1T}^*\right\|_2\\
\leq &\left\|B(\mathbf{X}_\mathcal{T})^\top w(\mathbf{X}_\mathcal{T}) - B(\mathbf{X}_\mathcal{T})^\top
B(\mathbf{X}_\mathcal{T})\lambda_1^{*}\right\|_2\left\|[B(\mathbf{X}_\mathcal{T})^\top
B(\mathbf{X}_\mathcal{T})]^{-1}\right\|_2.
\end{align*}

The first equality is due to the definition of $\lambda^\dagger$. The
second inequality is due to the Cauchy-Schwarz inequality.

We first consider the term $\left\|B(\mathbf{X}_\mathcal{T})^\top w(\mathbf{X}_\mathcal{T}) -
B(\mathbf{X}_\mathcal{T})^\top B(\mathbf{X}_\mathcal{T})\lambda_1^{*}\right\|_2$.
\begin{align*}
&\left\|B(\mathbf{X}_\mathcal{T})^\top w(\mathbf{X}_\mathcal{T}) - B(\mathbf{X}_\mathcal{T})^\top
B(\mathbf{X}_\mathcal{T})\lambda_1^{*}\right\|_2\\
\leq & \left\|B(\mathbf{X}_\mathcal{T})^\top w(\mathbf{X}_\mathcal{T}) - B(\mathbf{X}_\mathcal{T})^\top (1/\{n\pi(\mathbf{X}_\mathcal{T})\})\right\|_2 \\
&+ \left\|B(\mathbf{X}_\mathcal{T})^\top (1/\{n\pi(\mathbf{X}_\mathcal{T})\}) - B(\mathbf{X}_\mathcal{T})^\top B(\mathbf{X}_\mathcal{T})\lambda_1^{*}\right\|_2\\
=& \left\|B(\mathbf{X}_\mathcal{T})^\top w(\mathbf{X}_\mathcal{T}) - B(\mathbf{X}_\mathcal{T})^\top (1/\{n\pi(\mathbf{X}_\mathcal{T})\})\right\|_2 + \frac{1}{n}\left\|B(\mathbf{X}_\mathcal{T})\right\|_2\cdot O(K^{-r_\pi})\\
=& \left\|B(\mathbf{X}_\mathcal{T})^\top w(\mathbf{X}_\mathcal{T}) - \frac{1}{n} B(\mathbf{X})^\top 1\right\|_2 + \left\|\frac{1}{n} B(\mathbf{X})^\top 1 - B(\mathbf{X}_\mathcal{T})^\top (1/\{n\pi(\mathbf{X}_\mathcal{T})\})\right\|_2 \\
&+ \frac{1}{n}\left\|B(\mathbf{X}_\mathcal{T})\right\|_2\cdot O(K^{-r_\pi})\\
\leq & \left\|\delta\right\|_2\frac{\sum_{i=1}^n (1-Z_i)}{n} + \left\|\frac{1}{n} \sum_{i=1}^n B(X_i) - \frac{1}{n} \sum_{i=1}^n\frac{Z_i}{\pi(X_i)}B(X_i)\right\|_2 \\
&+ \frac{1}{n}\left\|B(\mathbf{X}_\mathcal{T})\right\|_2\cdot O(K^{-r_\pi})\\
\leq & \left\|\delta\right\|_2\frac{\sum_{i=1}^n (1-Z_i)}{n} + \left\|\frac{1}{n} \sum_{i=1}^n\frac{\pi(X_i)-Z_i}{\pi(X_i)}B(X_i)\right\|_2 + O(K^{1/2-r_\pi})\\
\leq & \left\|\delta\right\|_2 + \left\|\frac{1}{n} \sum_{i=1}^n\frac{\pi(X_i)-Z_i}{\pi(X_i)}B(X_i)\right\|_2 + O(K^{1/2-r_\pi}).
\end{align*}

We next bound the middle term $\left\|\frac{1}{n}
\sum_{i=1}^n\frac{\pi(X_i)-Z_i}{\pi(X_i)}B(X_i)\right\|_2$ by
Bernstein's inequality. Recall that the Bernstein's inequality for
random matrices \citep{tropp2015introduction} says, let $\{W_k\}$ be a
sequence of independent random matrices with dimensions $d_1\times
d_2$. Assume that $\mathbb{E}W_k = 0$ and $\left\|W_k\right\|_2\leq
R_n$ almost surely. Define
\[\sigma_n^2 = \max\{\left\|\sum^n_{k=1}\mathbb{E}(W_kW_k^\top)\right\|_2,
\left\|\sum^n_{k=1}\mathbb{E}(W_k^\top W_k)\right\|_2\}.\] Then for all $t\geq
0$,
\[P(\left\|\sum^n_{k=1}W_k\right\|_2\geq t)\leq
(d_1+d_2)\exp(-\frac{t^2/2}{\sigma^2_n+R_nt/3}).\]

Therefore, to bound $\left\|1 - \frac{1}{n} \sum_{i=1}^n\frac{Z_i}{\pi(X_i)}\right\|_2$,
we first show the summand is mean zero
\begin{align*}
\mathbb{E}\left\{\frac{1}{n} \sum_{i=1}^n\frac{\pi(X_i)-Z_i}{\pi(X_i)}B(X_i)\right\} = 0.
\end{align*} 

Furthermore, we have
\begin{align*} 
&\left\|\frac{1}{n} \sum_{i=1}^n\frac{\pi(X_i)-Z_i}{\pi(X_i)}B(X_i)\right\|_2\\
\leq &\frac{1}{n}\left\| \frac{\pi(X_i)-Z_i}{\pi(X_i)}\right\|_2\left\|B(X_i)\right\|_2\\
\leq &C'\frac{K^{1/2}}{n},
\end{align*} 
for some constant $C'$. The first inequality is due to Cauchy-Schwarz
inequality. The second equality is due to \Cref{assumption:propensity}.1.

Finally, we have
\begin{align*}
&\left\|\sum_{i=1}^n\mathbb{E}\{\frac{1}{n^2}(\frac{\pi(X_i)-Z_i}{\pi(X_i)})^2B(X_i)B(X_i)^\top\}\right\|_2\\
\leq &\frac{1}{n}\sup_i\left(\frac{\pi(X_i)-Z_i}{\pi(X_i)}\right)^2\left\|E\{B(X_i)B(X_i)^\top\}\right\|_2\\
\leq &\frac{C_1}{nC_3^2} = \frac{C''}{n},
\end{align*} 
for some constant $C''$. The last inequality is due to
\Cref{assumption:basis}.2 and \Cref{assumption:propensity}.1.

Therefore, by the Bernstein's inequality, we have
\begin{align*} P(\left\|
\frac{1}{n} \sum_{i=1}^n\frac{\pi(X_i)-Z_i}{\pi(X_i)}B(X_i)\right\|_2 >t)\leq
(K+1)\exp\left(-\frac{t^2/2}{\frac{C''}{n}+C'\frac{K^{1/2}}{n}\cdot t/3}\right).
\end{align*} 

When $t = O_p(K^{1/2}(\log K)/n)$, then the right side of the
inequality going to zero as $n\rightarrow\infty$ and for sufficiently
large constants $C'$ and $C''$. This gives
\begin{align*} \left\|
\frac{1}{n} \sum_{i=1}^n\frac{\pi(X_i)-Z_i}{\pi(X_i)}B(X_i)\right\|_2
= O_p(K^{1/2}(\log K)/n).
\end{align*}

Hence, we have
\begin{align*} 
&\left\|B(\mathbf{X}_\mathcal{T})^\top w(\mathbf{X}_\mathcal{T}) - B(\mathbf{X}_\mathcal{T})^\top
B(\mathbf{X}_\mathcal{T})\lambda_1^{*}\right\|_2\\
\leq & \left\|\delta\right\|_2 + \left\|\frac{1}{n} \sum_{i=1}^n\frac{\pi(X_i)-Z_i}{\pi(X_i)}B(X_i)\right\|_2 + O(K^{1/2-r_\pi})\\
\leq & \left\|\delta\right\|_2 + O_p(K^{1/2}(\log K)/n) +
O(K^{1/2-r_\pi})\\
\leq & \left\|\delta\right\|_2 + O_p(K^{1/2}((\log K)/n+K^{-r_\pi})\\
\leq &O_p(K^{1/2}((\log K)/n+K^{-r_\pi}).
\end{align*} 

Finally, this gives
\begin{align*} 
&\left\|[B(\mathbf{X}_\mathcal{T})^\top B(\mathbf{X}_\mathcal{T})]^{-1}[B(\mathbf{X}_\mathcal{T})^\top
w(\mathbf{X}_\mathcal{T})] - \lambda_{1T}^*\right\|_2\\ 
\leq &\left\|[B(\mathbf{X}_\mathcal{T})^\top w(\mathbf{X}_\mathcal{T})] - B(\mathbf{X}_\mathcal{T})^\top
B(\mathbf{X}_\mathcal{T})\lambda_1^{*}\right\|_2\left\|[B(\mathbf{X}_\mathcal{T})^\top
B(\mathbf{X}_\mathcal{T})]^{-1}\right\|_2\\
\leq &O_p(K^{1/2}((\log K)/n+K^{-r_\pi})\cdot 2/(C_2n)\\ 
=&O_p(K^{1/2}((\log K)/n+K^{-r_\pi})/n),\label{eq:suppcorrrate}
\end{align*} 
where the second inequality is due to the following calculation
\begin{align*} 
&\lambda_{\min}(B(\mathbf{X}_\mathcal{T})^\top
B(\mathbf{X}_\mathcal{T})/\sum_{i=1}^nZ_i)\\
\geq &\lambda_{\min}(\mathbb{E}[B(X_i)^\top
B(X_i)]) - \left\|B(\mathbf{X}_\mathcal{T})^\top
B(\mathbf{X}_\mathcal{T})/\sum_{i=1}^nZ_i-\mathbb{E}[B(X_i)^\top
B(X_i)]\right\|_2\\
\geq & C_2 - \left\|B(\mathbf{X}_\mathcal{T})^\top
B(\mathbf{X}_\mathcal{T})/\sum_{i=1}^nZ_i-\mathbb{E}[B(X_i)^\top
B(X_i)]\right\|_2\\
\geq & C_2/2
\end{align*}
for sufficiently large $n$. The first inequality is due to the Weyl
inequality. The second inequality is due to
\Cref{assumption:basis}. The third inequality is due to
$\left\|B(\mathbf{X}_\mathcal{T})^\top B(\mathbf{X}_\mathcal{T})/\sum_{i=1}^nZ_i-\mathbb{E}[B(X_i)^\top
B(X_i)]\right\|_2 = o_p(1)$ from Lemma D.4 of \citet{fan2016improving}.

To conclude the proof, we have
\begin{align*} 
&sup_{x\in\mathcal{X}}\left|nw(x) - \frac{1}{\pi(x)}\right|
\\
\leq&sup_{x\in\mathcal{X}}|nw(x) - nB(x)\lambda^{\dagger}| + \sup_{x\in\mathcal{X}}|nB(x)\lambda^{\dagger}-nB(x)\lambda_{1T}^*| + |nB(x)\lambda_{1T}^* - \frac{1}{\pi(x)}|\\
=& O_p(K((\log K)/n+K^{-r_\pi}))+O_p(K^{1/2}((\log K)/n+K^{-r_\pi})/n)\cdot O_p(K^{1/2}n)\\
&+ o_p(K^{-r_\pi})\\
=&O_p(K((\log K)/n+K^{-r_\pi}))+ o_p(K^{-r_\pi})\\ 
=&O_p(K((\log K)/n+K^{-r_\pi}))\\ 
=&o_p(1).
\end{align*} 

The first inequality is due to the triangle inequality. The second
equality is due to the piece-wise constant condition on $w(\cdot)$ and
\Cref{assumption:basis}.1. The last equality is due to
\Cref{assumption:hyperparam}.1.

Similarly, we have
\begin{align*} 
&\left\|nw(x) - \frac{1}{1-\pi(x)}\right\|_2 \\
\leq & O_p(K^{1/2}((\log K)/n+K^{-r_\pi})/n)\cdot O_p(K^{1/2}n)+ o_p(K^{-r_\pi})\\
=&O_p(K((\log K)/n+K^{-r_\pi}))\\ 
=&o_p(1).
\end{align*}
\hfill
\QED

Analogously, we can prove a similar lemma for the weights of the
control units.

\begin{lemma}
\label{lemma:wtconsistencycontrol} Let $\pi(x)$ be the true propensity score
function. Under \Cref{assumption:basis} and
\Cref{assumption:propensity}, we have, for all control units,
\begin{enumerate}
    \item $\sup_{x\in\mathcal{X}}|nw(x) - 1/(1-\pi(x))| =O_p(K((\log K)/n+K^{-r_\pi}))=o_p(1)$,
    \item $\left\|nw(x) - 1/(1-\pi(x))\right\|_{P,2} =O_p(K((\log K)/n+K^{-r_\pi}))= o_p(1)$.
\end{enumerate}
\end{lemma}

Building on \Cref{lemma:wtconsistency} and
\Cref{lemma:wtconsistencycontrol}, we then establish the asymptotic
normality and semiparametric efficiency for the \gls{ATE} estimator
\Cref{thm:normality}.

\paragraph{Proof.}

The proof utilizes empirical processes techniques as in
\citet{fan2016improving}.

We first decompose $\hat{\mu} - \mu$ into a main term and a few
residual terms:
\begin{align*}
&\hat{\mu} - \mu \\
= & \sum_{i=1}^nZ_iw(X_i)Y_i - \sum_{i=1}^n(1-Z_i)w(X_i)Y_i - (E\{Y_i(1)\} - E\{Y_i(0)\})\\
= & \left(\sum_{i=1}^nZ_iw(X_i)Y_i - E\{Y_i(1)\}\right) - \left(\sum_{i=1}^n(1-Z_i)w(X_i)Y_i - E\{Y_i(0)\}\right) \\
=&\frac{1}{n}\sum_{i=1}^n\left[\frac{Z_i}{\pi(X_i)}(Y_i -Y_1(X_i)) - \frac{1-Z_i}{1-\pi(X_i)}(Y_i-Y_0(X_i)) +  (Y_1(X_i)-Y_0(X_i)) - \mu\right]\\
&+\sum_{i=1}^n(w(X_i)-\frac{1}{n(1-\pi(X_i))})(1-Z_i)(Y_i-Y_0(X_i))\\
&+\sum_{i=1}^n(w(X_i)-\frac{1}{n\pi(X_i)})Z_i(Y_i-Y_1(X_i))\\
&+\sum_{i=1}^n(w(X_i)(1-Z_i)-\frac{1}{n})(Y_0(X_i)-B(X_i)^\top\lambda_{C}^* )\\
&+\sum_{i=1}^n(w(X_i)Z_i- \frac{1}{n})(Y_1(X_i)-B(X_i)^\top\lambda_{T}^* )\\
& + \sum_{i=1}^n(w(X_i)(1-Z_i)-\frac{1}{n})(B(X_i)^\top\lambda_{2C}^* )\\
& + \sum_{i=1}^n(w(X_i)Z_i-\frac{1}{n})(B(X_i)^\top\lambda_{2T}^* )\\
= &\frac{1}{n} \sum^n_{i=1}S_i + R_{0C} + R_{0T} + R_{1C} + R_{1T} + R_{2C} + R_{2T},
\end{align*}
where
\begin{align*}
S_i = &\frac{Z_i}{\pi(X_i)}(Y_i -Y_1(X_i)) - \frac{1-Z_i}{1-\pi(X_i)}(Y_i-Y_0(X_i)) +  (Y_1(X_i)-Y_0(X_i)) - \mu\\
R_{0C} = &\sum_{i=1}^n(w(X_i)-\frac{1}{n(1-\pi(X_i))})(1-Z_i)(Y_i-Y_0(X_i)),\\
R_{0T} = &\sum_{i=1}^n(w(X_i)-\frac{1}{n\pi(X_i)})Z_i(Y_i-Y_1(X_i)),\\
R_{1C} = &\sum_{i=1}^n(w(X_i)(1-Z_i)-\frac{1}{n})(Y_0(X_i)-B(X_i)^\top\lambda_{C}^* ), \\
R_{1T} = &\sum_{i=1}^n(w(X_i)Z_i- \frac{1}{n})(Y_1(X_i)-B(X_i)^\top\lambda_{T}^* ),\\
R_{2C} = &\sum_{i=1}^n(w(X_i)(1-Z_i)-\frac{1}{n})(B(X_i)^\top\lambda_{2C}^* ),\\
R_{2T} = &\sum_{i=1}^n(w(X_i)Z_i-\frac{1}{n})(B(X_i)^\top\lambda_{2T}^* ).
\end{align*}

Below we show $R_{jT} = o_p(n^{-1/2}), R_{jC} = o_p(n^{-1/2}), 0 \leq
j \leq 2$. The conclusion follows from $S_i$ taking the same form as
the efficient score \citep{hahn1998role}. $\hat{\mu}$ is thus
asymptotically normal and semi-parametrically efficient. 

Given \Cref{lemma:wtconsistency} and
\Cref{lemma:wtconsistencycontrol}, the rest of the proof repeats the
proof structure of Theorem 3 in \citet{wang2020minimal}. They prove a
similar consistency and asymptotic normality result for weighting
methods. We leverage \Cref{lemma:wtconsistency} and
\Cref{lemma:wtconsistencycontrol} to adapt their proof to matching for
balance. We include the rest of the proof here for completeness.

The first term we study is
$R_{0C}=\sum_{i=1}^n(w(X_i)-\frac{1}{n(1-\pi(X_i))})(1-Z_i)(Y_i-Y_0(X_i))$.
Consider an empirical process $\mathbb{G}_n(f_0) = n^{1/2}
(\mathbb{P}_n - \mathbb{P}) f_0(Z, Y, X)$, where $\mathbb{P}_n$ stands
for the empirical measure and $\mathbb{P}$ stands for the expectation,
and
\[f_0(Z, Y, X) = (nw(X)-\frac{1}{1-\pi(X)})(1-Z)(Y-Y_0(X)).\]
By the unconfoundedness assumption, we have that $\mathbb{P}f_0(Z,
Y,X) = 0.$ By Markov's inequality and maximal inequality, we have
\[\sqrt{n}R_0\leq\sup_{f_0\in\mathcal{F}}\mathbb{G}_n(f_0)\lesssim
\mathbb{E}\sup_{f_0\in\mathcal{F}}\mathbb{G}_n(f_0)\lesssim
J_{[]}(\left\|F_0\right\|_{P,2},\mathcal{F},L_2(P)),\] where the set
of functions is
$\mathcal{F}=\{f_0:\left\|w(\cdot)-1/n(1-\pi(\cdot))\right\|_\infty\leq\delta\}$,
where $\left\|f\right\|_{\infty}=\sup_{x\in\mathcal{X}}|f(x)|$ and
$\delta=C(K((\log K)/n+K^{-r_\pi}))$ for some constant $C>0$. The
second inequality is due to Markov's inequality.
$J_{[]}(\left\|F_0\right\|_{P,2},\mathcal{F},L_2(P)$ is the bracketing
integral. $F_0:= \delta|Y-Y_0(X)|\gtrsim |f_0(Z, Y, X)|$ is the
envelop function. We also have
$\left\|F_0\right\|_{P,2}=(\mathbb{P}F_0^2)^{1/2}\lesssim\delta$ by
$\mathbb{E}|Y-Y_0(X)|<\infty.$ This construction is due to
\Cref{lemma:wtconsistency}, where we have
\[\left\|nw(x) - \frac{1}{1-\pi(x)}\right\|_2 \\
=O_p(K((\log K)/n+K^{-r_\pi}))=o_p(1).\]

We next show $\mathbb{E}\sup_{f_0\in\mathcal{F}}\mathbb{G}_n(f_0)
\rightarrow 0$ as $\delta\rightarrow 0$. We bound $J_{[]}(\left\|F_0\right\|_{P,2},\mathcal{F},L_2(P)$ by $n_{[]}(\varepsilon,\mathcal{F}, L_2(P))$:
\[J_{[]}(\left\|F_0\right\|_{P,2},\mathcal{F},L_2(P) \lesssim \int^\delta_0 \sqrt{n_{[]}(\varepsilon,\mathcal{F}, L_2(P))}\text{d}\varepsilon.\]

Then we have
\begin{align*}
\log n_{[]}(\varepsilon,\mathcal{F}, L_2(P))&\lesssim \log n_{[]}(\varepsilon,\mathcal{F}_0\delta, L_2(P))\\
&=\log n_{[]}(\varepsilon/\delta,\mathcal{F}_0, L_2(P))\\
&\lesssim \log n_{[]}(\varepsilon/\delta,\mathcal{M}, L_2(P)\\
&\lesssim(\delta/\varepsilon)^{(1/k_1)},
\end{align*}
where we consider a new set of functions
$\mathcal{F}_0=\{f_0:\left\|w(\cdot)-1/n(1-\pi(\cdot))\right\|_\infty\leq
C\}$ for some constant $C>0$. The first inequality is due to
$w(\cdot)$ bounded away from 0 and (can always be made) Lipschitz. The
last inequality is due to \Cref{assumption:outcome}.4.

Therefore, we have
\[J_{[]}(\left\|F_0\right\|_{P,2},\mathcal{F},L_2(P)\lesssim\int^\delta_0\sqrt{\log
n_{[]}(\varepsilon,\mathcal{F},
L_2(P))}d\varepsilon\lesssim\int^\delta_0(\delta/\varepsilon)^{(1/2k_1)}d\varepsilon.\]

This goes to 0 as $\delta$ goes to 0 by $2k_1>1$ and the integral
converges. Thus, this shows that $n^{1/2}R_{0C}=o_p(1)$. With the
exact same argument, we can also show that $n^{1/2}R_{0T}=o_p(1)$.\\

Next, we consider $R_{1C}
=\sum_{i=1}^n(w(X_i)(1-Z_i)-\frac{1}{n})(Y_0(X_i)-B(X_i)^\top\lambda_{C}^*
)$. Define the empirical process
$\mathbb{G}_n(f_1)=n^{1/2}(\mathbb{P}_n-\mathbb{P})f_1(Z,X)$, where
\[f_1(Z,X)=(n(1-Z)w(X)-1)(Y_0(X)-B(X)^\top\lambda_{2C}^*).\]
We have
\begin{align*}
n^{1/2}R_1 &= \mathbb{G}_n(f_1)+n^{1/2} \mathbb{P}f_1(Z,X)\\
&\leq\sup_{f_1\in\mathcal{F}_1}\mathbb{G}_n(f_1)+n^{1/2}\sup_{f_1\in\mathcal{F}_1}\mathbb{P}f_1,
\end{align*}
where $\Delta(X) := Y_0(X_i)-B(X_i)^\top\lambda_{C}^*$, $\mathcal{F}_1
= \{f_1:\left\|w(\cdot)-1/n(1-\pi(\cdot))\right\|_{P,2}\leq \delta_1,
\left\|\Delta\right\|_\infty\leq\delta_2\}, \delta_1 = C
(K((\log K)/n+K^{-r_\pi}))$, and $\delta_2 = CK^{-r_y}$ for some
constant $C > 0.$ This construction is due to
\Cref{lemma:wtconsistencycontrol}.2 and \Cref{assumption:outcome}.3.

Similar to characterizing $R_{0C}$, we have
\begin{align*}
\sup_{f_1\in\mathcal{F}_1}\mathbb{G}_n(f_1)\lesssim \mathbb{E}\sup_{f_1\in\mathcal{F}_1}\mathbb{G}_n(f_1)\lesssim J_{[]}(\left\|F_1\right\|_{P,2},\mathcal{F},L_2(P)),
\end{align*}
where $F_1:=C\delta_2$ for some constant $C > 0$ so that
$\left\|F_1\right\|_{P,2} \lesssim \delta_2.$ We then bound
$J_{[]}(\left\|F_1\right\|_{P,2},\mathcal{F}_1,L_2(P))$:
$J_{[]}(\left\|F_1\right\|_{P,2},\mathcal{F}_1,L_2(P) \lesssim \int^\delta_0 \sqrt{n_{[]}(\varepsilon,\mathcal{F}_1, L_2(P))}\text{d}\varepsilon.$ Moreover, we have
\begin{align*}
\log n_{[]}(\varepsilon, \mathcal{F}_1, L_2(P))&\lesssim \log n_{[]}(\varepsilon/\delta_2, \mathcal{F}_0, L_2(P))\\
&\lesssim \log n_{[]}(\varepsilon/\delta_1, G_{10}, L_2(P)) + \log n_{[]}(\varepsilon/\delta_2, G_{20}, L_2(P))\\
&\lesssim \log n_{[]}(\varepsilon/\delta_1, \mathcal{M}, L_2(P)) + \log n_{[]}(\varepsilon/\delta_2, \mathcal{H}, L_2(P))\\
&\lesssim (\delta_1/\varepsilon)^{1/k_1}+(\delta_2/\varepsilon)^{1/k_2},
\end{align*}
where
$\mathcal{F}_0=\{f_1:\left\|w(\cdot)-1/(1-\pi(\cdot))\right\|_{P,2}\leq
C, \left\|\Delta\right\|_{P,2}\leq 1\}, \mathcal{G}_{10}=\{m\in
\mathcal{M}+1/(1-\pi(\cdot)):\left\|w\right\|_{P,2}\leq
C\},\mathcal{G}_{20}=\{\Delta\in\mathcal{H}-B(X_i)^\top\lambda_{2C}^*:\left\|\Delta\right\|_{P,2}\leq
1\}.$ 

Therefore we have
\[J_{[]}(\left\|F_1\right\|_{P,2},\mathcal{F}_1,L_2(P)\lesssim\int^\delta_0(\delta_1/\varepsilon)^{(1/2k_1)}d\varepsilon + \int^\delta_0(\delta_2/\varepsilon)^{(1/2k_2)}d\varepsilon.\]
By $2k_1>1, 2k_2>1$ (\Cref{assumption:outcome}.4), we have
$J_{[]}(\left\|f_1\right\|_{P,2}, \mathcal{F}, L_2(P)) = o(1)$. This
gives $\sup_{f_1\in\mathcal{F}_1}\mathbb{G}_n(f_1) = o_p(1)$.

Finally we show that $n^{1/2}\sup_{f_1\in\mathcal{F}_1}\mathbb{P}f_1 =
o_p(1).$
\begin{align*}
n^{1/2}\sup_{f_1\in\mathcal{F}} \mathbb{P}f_1 
&=n^{1/2}\sup_{m\in\mathcal{G}_1, \Delta\in\mathcal{G}_2}\mathbb{E}(n(1-Z)(w(X)-Z/\pi(X))\Delta(X))\\
&=n^{1/2}\sup_{m\in\mathcal{G}_1, \Delta\in\mathcal{G}_2}\mathbb{E}((nw(X)(1-\pi(X))-1)\Delta(x))\\
&\lesssim n^{1/2} \sup_{m\in\mathcal{G}_1}\left\|nw(x)-\frac{1}{1-\pi(x)}\right\|_{P,2} \sup_{\Delta\in\mathcal{G}_2}\left\|\Delta(x)\right\|_{P,2}\\
&\lesssim n^{1/2}\delta_1\delta_2 = o_p(1),
\end{align*}
where $\mathcal{G}_1=\{m\in\mathcal{M}:
\left\|w(\cdot)-1/(1-\pi(\cdot))\right\|_{P,2}\leq\delta_1\},
\mathcal{G}_2=\{\Delta\in\mathcal{H}-B(X_i)^\top\lambda_{2C}^*:\left\|\Delta\right\|_\infty\leq\delta_2\}$.
The last equality is due to assumption $n^{1/2}\lesssim
K^{r_\pi+r_y-1/2}$.

Therefore, we can conclude $n^{1/2}R_{1C}=o_p(1).$ Analogously, we can
conclude $n^{1/2}R_{1T}=o_p(1).$

Lastly, we study $R_{2C} =
\sum_{i=1}^n(w(X_i)(1-Z_i)-\frac{1}{n})(B(X_i)^\top\lambda_{2C}^* )$. We
have $||R_{2C}||_2^2 = o_p(n^{-1/2})$ because the covariate balancing
condition in the optimization problem and
\Cref{assumption:hyperparam}.2. Similarly, we have $R_{2T} =
o_p(n^{-1/2}).$

After establishing consistency and asymptotic normality, we provide an
estimator for the asymptotic variance. We construct the variance
estimator by approximating the efficient influence function. This
construction and the proof mimic the second part of Theorem 3 in
\citet{wang2020minimal}.

First recall that the semiparametric efficiency bound for \gls{ATE}
\citep{hahn1998role} is
\begin{align}
V_{opt}:= &\mathbb{E}\left[\frac{var(Y_i(1)\mid
X_i))}{\pi(X_i)} + \frac{var(Y_i(0)\mid
X_i))}{1-\pi(X_i)} +(E(Y_i(1)-Y_i(0)\mid
X_i)-\mu)^2\right]\\
=&\mathbb{E}\left[\left(\frac{Z_iY_i}{\pi(X_i)} - \bar{Y}_{i}(1)\right) - \left(\frac{(1-Z_i)Y_i}{1-\pi(X_i)} - \bar{Y}_{i}(0)\right) \right.\nonumber\\
&\left.- \left(Y_1(X_i)(\frac{Z_i}{\pi(X_i)}-1) - Y_0(X_i)(\frac{1-Z_i}{1-\pi(X_i)} - 1)\right) \right]^2
\end{align}

To estimate this variance, we consider the following estimator
\begin{align*}
\begin{split}
\hat{V}_K = &\sum^n_{i=1} \left[ Z_iw(X_i)Y_i -
\frac{\sum_{i=1}^nZ_iw(X_i)Y_i}{\sum_{i=1}^nZ_i} \right.\\
&- (1-Z_i)w(X_i)Y_i +
\frac{\sum_{i=1}^n(1-Z_i)w(X_i)Y_i}{\sum_{i=1}^n(1-Z_i)} \\
& - \hat{Y}_T(X_i)\times \left(Z_iw(X_i)
- \frac{1}{\sum_{i=1}^nZ_i}\right)\\
& \left.+ \hat{Y}_C(X_i)\times \left((1-Z_i)w(X_i)
- \frac{1}{\sum_{i=1}^n(1-Z_i)}\right) \right]^2.
\end{split}
\end{align*} 
where
\begin{multline*}
\hat{Y}_T(X_i) = B(X_i)^\top
\left\{\frac{\sum_{i=1}^nZ_iw(X_i)B(X_i) B(X_i)^\top}{\sum_{i=1}^nZ_i}
\right\}^{-1}
\cdot\left\{ \frac{\sum_{i=1}^nZ_iw(X_i)B(X_i) Y_i}{\sum_{i=1}^nZ_i} \right\}
\end{multline*}
and
\begin{multline*}
\hat{Y}_C(X_i)  = B(X_i)^\top \left\{\frac{\sum_{i=1}^n(1-Z_i)w(X_i)B(X_i) B(X_i)^\top}{\sum_{i=1}^n(1-Z_i)}
\right\}^{-1}
\cdot\left\{ \frac{\sum_{i=1}^n(1-Z_i)w(X_i)B(X_i) Y_i}{\sum_{i=1}^n(1-Z_i)} \right\}
\end{multline*}
are least square estimators of $Y_1(X_i)$ and $Y_0(X_i)$ respectively.

To prove consistency, it is sufficient show that 
\begin{align}
\hat{Y}_T(X_i) - Y_1(X_i)&\stackrel{a.s.}{\rightarrow} 0,\label{eq:lsq-consist}\\
\hat{Y}_C(X_i) - Y_0(X_i)&\stackrel{a.s.}{\rightarrow} 0,\label{eq:lsq-consist2}
\end{align}
because we have shown that $nw(X_i)$ is consistent for $\pi(X_i)$ and
$\frac{\sum_{i=1}^nZ_iw(X_i)Y_i}{\sum_{i=1}^nZ_i} -
\frac{\sum_{i=1}^n(1-Z_i)w(X_i)Y_i}{\sum_{i=1}^n(1-Z_i)}$ is
consistent for $\mu$.

Below we assume a stronger smoothness assumption, i.e. $r_y>1$.

To prove \Cref{eq:lsq-consist}, we first rewrite $Y_i$ as $Y_1(X_i) =
B(X_i)^\top\lambda_{2T}^* + \gamma + \epsilon_i,$ where $\gamma =
O(K^{-r_y})$ from \Cref{assumption:outcome}.3, and $\epsilon_i$ is
some iid zero mean error with variance $\sigma^2 = var(Y_1(X_i)|X_i)$.
Therefore,
\begin{align*}
&\left\{\frac{\sum_{i=1}^nZ_iw(X_i)B(X_i) B(X_i)^\top}{\sum_{i=1}^nZ_i}
\right\}^{-1}\cdot\left\{ \frac{\sum_{i=1}^nZ_iw(X_i)B(X_i) Y_i}{\sum_{i=1}^nZ_i} \right\}\\
=&\left\{\frac{\sum_{i=1}^nZ_iw(X_i)B(X_i) B(X_i)^\top}{\sum_{i=1}^nZ_i}
\right\}^{-1}\cdot\left\{ \frac{\sum_{i=1}^nZ_iw(X_i)B(X_i)^\top \{B(X_i)^\top\lambda_2^* + \gamma + \epsilon_i\}}{\sum_{i=1}^nZ_i} \right\}
\\
=&\lambda^*_2 + \left\{\frac{\sum_{i=1}^nZ_iw(X_i)B(X_i) B(X_i)^\top}{\sum_{i=1}^nZ_i}
\right\}^{-1}\cdot\left\{ \frac{\sum_{i=1}^nZ_iw(X_i)B(X_i)^\top \{B(X_i)^\top\gamma\}}{\sum_{i=1}^nZ_i}\right\}\\
&+\left\{\frac{\sum_{i=1}^nZ_iw(X_i)B(X_i) B(X_i)^\top}{\sum_{i=1}^nZ_i}
\right\}^{-1}\cdot\left\{ \frac{\sum_{i=1}^nZ_iw(X_i)B(X_i)^\top \epsilon_i}{\sum_{i=1}^nZ_i} \right\}\\
=&\lambda^*_2 + E\{Z_iw(X_i)B(X_i)^\top
B(X_i)\}^{-1}E\{Z_iw(X_i)B(X_i)^\top\}\{\gamma+E(\epsilon_i)\} + O_p(n^{-1/2})\\
=&\lambda^*_{2T} + O_p(K^{-r_{y}+1/2}),
\end{align*}
where the last equality is due to \Cref{assumption:basis}.4 and \Cref{assumption:outcome}.3 and the law of large numbers.

Finally we have 
\begin{align*}
&B(X_i)^\top\left\{\frac{\sum_{i=1}^nZ_iw(X_i)B(X_i) B(X_i)^\top}{\sum_{i=1}^nZ_i}
\right\}^{-1}\cdot\left\{ \frac{\sum_{i=1}^nZ_iw(X_i)B(X_i) Y_i}{\sum_{i=1}^nZ_i} \right\}\\
=&B(X_i)^\top\lambda^*_{2T} + B(X_i) \cdot O_p(K^{-r_{y}+1/2})\\
=&Y_1(X_i) + B(X_i) \cdot O_p(K^{-r_{y}+1/2}) + O_p(K^{-r_{y}})\\
=&Y_1(X_i) + o_p(1)
\end{align*}
The last equality is due to assumption
\Cref{assumption:basis}.4 and the additional assumption $r_y>1$.

\Cref{eq:lsq-consist2} can be proved analogously. Hence we
established the consistency of $\hat{V}_K$.

\glsresetall

\glsresetall

\section{Existence of a solution}

In this section, we provide sufficient conditions that guarantee the existence of
a solution to the matching for balance optimization problem. 
As discussed, if there is no solution for the matching for balance problem, there will likely be no solution that balances covariates also for other matching or weighting methods with positive weights, i.e. that do not extrapolate.
Together with \Cref{thm:normality}, this result describes the settings where matching methods are as statistically efficient as weighting methods. 
Moreover, they are both semiparametrically efficient.

\begin{proposition}(Sufficient conditions for the existence of a
solution to matching for balance)
\label{prop:solutionexist}
If there exists a constant $C_3 > 0$ such that $C_3 < \pi(x) < 1-C_3$
and
\begin{align}
\label{eq:assumption}
C_3\geq \Theta(1 / ((\log K)+nK^{-r_\pi})),
\end{align}
then matching for balance admits a solution with probability
$1-\delta_0$ under \Cref{assumption:basis,assumption:propensity,assumption:outcome,assumption:hyperparam} when $n
\geq \log_{1-\rho}(\delta_0  2^{-K})$ for a constant
$\rho \in (0,1)$. (The precise technical definition of the constant $\rho$ is in
\Cref{sec:solutionexistproof} of the Supplementary Materials.)
\end{proposition}

\Cref{prop:solutionexist} describes sufficient conditions for a
matching-for-balance solution to exist. It roughly requires that the
propensity score function should be bounded away from zero and one; it
must be at least $\Theta(1 / ((\log K)+nK^{-r_\pi}))$ away. As the
number of units $n$ increases, this requirement becomes increasingly
weak and a match-for-balance solution exists more likely. Together
with \Cref{thm:normality}, \Cref{prop:solutionexist} delineates a
setting where matching methods can be as efficient as weighting
methods despite its integer constraints. In this setting, matching for
balance is asymptotically optimal: it is $\sqrt{n}$-consistent,
asymptotically normal, and semiparametrically efficient.

The intuition behind \Cref{prop:solutionexist} is twofold: (1) the
treated and control population are closer if the propensity score is
farther from zero and one; (2) if the two populations are closer, it
is more likely for a matching solution to exist, i.e. satisfy the
covariate balancing constraints in
\Cref{eq:matchcovbalcondition,eq:matchcovbalcondition_2}. Due to this
intuition, we posit the overlap condition in \Cref{prop:solutionexist}
to constrain how far away the two populations can be. It requires the
minimum propensity score $C_3$ to be larger than $\Theta(1 / ((\log
K)+nK^{-r_\pi}))$. This condition is stronger than what is usually
required of the overlap between the treated and the control (e.g.,
\Cref{assumption:propensity}\textcolor{black}{.1}). Nevertheless, it
guarantees the existence and asymptotic optimality of a matching
solution. The full proof of \Cref{prop:solutionexist} is in
\Cref{sec:solutionexistproof} of the Supplementary Material.

\section{Proof of Proposition 1}

\label{sec:solutionexistproof}

We first define the constant $\rho$ in \Cref{prop:solutionexist}. The
constant $\rho$ is defined as $\rho \stackrel{\Delta}{=}
\min_i P(X\in R_i\,|\, T=0) > 0$, where $R_i, i=1, \ldots, 3^K,$ are
the $3^K$ boxes centered at $E[B(X_i)\,|\,Z_i=1] +
\frac{3}{2}\delta\odot b$. The constant $\delta = (\delta_1, \ldots,
\delta_k)$ is the covariate imbalances allowed in matching for balance
(\Cref{eq:matchcovbalcondition,eq:matchcovbalcondition_2}) and
$b\in\mathbb{R}^K$ is a vector that each entry can be $-1, 0, 1$.

We will prove the existence of $w_C(X_j)$ that solves the optimization
problem of matching for balance, where $w_C(X_j) = (\sum
_{i=1}^nZ_im_{ij})/M.$ The exact same argument can establish the
existence of $w_T(X_j)$. 

The proof proceeds in three steps: (1) Show that a set of $w_C(X_i)$
exists that satisfies \Cref{eq:matchcovbalcondition}, without
conforming to the form of \Cref{eq:weights}. (2) Show that a set of
$w_C(X_i)$ exists that satisfies both \Cref{eq:matchcovbalcondition}
and \Cref{eq:weights}. (3) Show that a set of $m_{ij}$ satisfies both
\Cref{eq:matchcovbalcondition} and \Cref{eq:weights}.

The first step follows directly from Lemma 2 of
\citet{zhao2017entropy}. They show that a solution exists for
\Cref{eq:matchcovbalcondition} with probability $\delta_0$ when the
number of units satisfy $n \geq \log_{1-\rho}(\delta_0\cdot 2^{-K})$.
We call this solution $w_C^0(X_i)$.

The second step relies on the following observation: \Cref{eq:weights}
amounts to the restriction that  we can only approximate $w_C^0(X_i)$
up to the precision $(M)^{-1}$. The reason is that $w_C(X_i) = (\sum
_{i=1}^nZ_im_{ij})/M \in
\{k/M, k\in \mathbb{N}\}$; matching restricts both the numerator
$(\sum _{i=1}^nZ_im_{ij})$ and the denominator $M$ of $w_C(X_i)$ to be
integers.

However, we show that the $w_C(X_i) = (\sum _{i=1}^nZ_im_{ij})/M$
closest to $w_C^0(X_i)$ still satisfies
\Cref{eq:matchcovbalcondition}. Notice that $|w_C(X_i) - w_C^0(X_i)| <
1/M$. This implies
\begin{align*}
&\frac{1}{\sum_{i=1}^nZ_i}\left|\sum_{i=1}^nZ_iB_k(X_i) -
\sum_{j=1}^nw_C(X_j)(1-Z_j)B_k(X_j) \right|\\
\leq&\frac{1}{\sum_{i=1}^nZ_i}\left|\sum_{i=1}^nZ_iB_k(X_i) -
\sum_{j=1}^nw_C^0(X_j)(1-Z_j)B_k(X_j) \right| \\
&+\frac{1}{\sum_{i=1}^nZ_i}\left|\sum_{j=1}^nw_C^0(X_j)(1-Z_j)B_k(X_j)-
\sum_{j=1}^nw_C(X_j)(1-Z_j)B_k(X_j) \right|\\
\leq&\delta_k +
\frac{\sum_{i=1}^n(1-Z_i)}{\sum_{i=1}^nZ_i}\frac{1}{M}\left|\sum_{j=1}^nB_k(X_i)\right|\\
\stackrel{\Delta}{=}&\delta_k^M
\end{align*}
\Cref{assumption:hyperparam}.2 of \Cref{thm:normality} requires that
$||\delta||_2 = O_p(K^{1/2}((\log K)/n+K^{-r_\pi})$, where $\delta =
(\delta_1,\ldots, \delta_K).$ We will show that $||\delta^M||_2$ also
satisfies \Cref{assumption:hyperparam}.2. Therefore, $w_C(X_j)$ can
also lead to consistent, asymptotically normal, and semiparametrically
efficient \gls{ATE} estimators.

We bound the norm of $||\delta^M||_2$:
\begin{align*}
||\delta^M||_2 &\leq ||\delta||_2 + \sum_{k=1}^K \left[\frac{\sum_{i=1}^n(1-Z_i)}{\sum_{i=1}^nZ_i}\frac{1}{M}\left|\sum_{j=1}^nB_k(X_i)\right|\right]\\
&\leq ||\delta||_2 + \frac{C}{M}\frac{\sum_{i=1}^n(1-Z_i)}{\sum_{i=1}^nZ_i}K^{1/2} \\
&= O_p(K^{1/2}((\log K)/n+K^{-r_\pi}).
\end{align*}
The last equality is because we can let $M = \Theta(((\log
K)/n+K^{-r_\pi})^{-1})$. It is feasible because the number of match
required $M$ is smaller than the number of treated units or control
units:
\begin{align*}
M=&\Theta(((\log
K)/n+K^{-r_\pi})^{-1})\\
\leq &nC_3 - \sqrt{\log(1-\rho) + \frac{K\log 2}{n}}\\
\leq &nC_3 - \sqrt{(-\log \delta_0) / (2n)} \\
\leq &\min(\sum_{i=1}^nZ_i, \sum_{i=1}^n1-Z_i)
\end{align*}
with probability $1-\delta_0$. The first inequality is due to
\Cref{eq:assumption}. The second inequality is due to the assumption
that $n \geq
\log_{1-\rho}(\delta_0\cdot 2^{-K})$. The third inequality is due to
the Hoeffding's inequality.

The third step is to construct a solution $\{m_{ij}\}_{i,j=1}^n$ such
that they are consistent with the solution from the second step:
$w_C(X_j) = (\sum _{i=1}^nZ_im_{ij})/M$. We first notice that we only
need to make sure each treated unit is matched $M$ times. Moreover,
for each control unit $j$, it is matched $w_C(X_j)\cdot M$ times.
Therefore, to construct a solution $m_{ij}$, we match each control
unit to the first $w_C(X_j)\cdot M$ treated units that have not been
matched for $M$ times. More precisely, we start with setting
$m_{ij}=0$ for all $i,j$. We then iterate through the set of control
units $\{i:Z_i=0\}$. For each $i$, we set $m_{ij} = 1$ if $j$
satisfies $\sum_{k=1}^j
\mathbb{I}\{\sum_{j=1}^n m_{kj} < M\}< M$. Iterating through all the
control units leads to a solution $m_{ij}$ that satisfies both
\Cref{eq:matchcovbalcondition} and \Cref{eq:weights}.
\hfill\QED

\putbib[mybibliography17]
\end{bibunit}

\end{document}